\documentclass[12pt]{article}
\usepackage{amsmath,amssymb}

\newcommand{\ncm}{\newcommand}
 \ncm{\R}{\mathbb{R}}
 \ncm{\C}{\mathbb{C}}
 \ncm{\Q}{\mathbb{Q}}
 \ncm{\N}{\mathbb{N}}
 \ncm{\Ad}{\mbox{\rm Ad}}
 \ncm{\ad}{{\rm ad}}
 \ncm{\Ran}{{\rm Ran}}
 \ncm{\Aut}{\mbox{\rm Aut}}
 \ncm{\Sp}{{\rm Sp}}
 \ncm{\Spec}{{\rm Spec}}
 \ncm{\supp}{{\rm supp}}
 \ncm{\Ker}{{\rm Ker}}
 \ncm{\id}{{\rm id}}
 \ncm{\ra}{\rightarrow}
 \ncm{\cstar}{C$^*$-algebra}
 \ncm{\M}{{\mathcal{M}}}
 \ncm{\Ex}{{\mathcal Ex}}
 \ncm{\E}{{\mathcal E}}
 \ncm{\G}{{\mathcal G}}
 \ncm{\F}{{\mathcal F}}
 \ncm{\B}{{\mathcal B}}
 \ncm{\T}{\mathbb{T}}
 \ncm{\V}{{\mathcal V}}
 \ncm{\K}{{\mathcal K}}
 \ncm{\CO}{{\mathcal O}}
 \ncm{\D}{{\mathcal D}}
 \ncm{\I}{{\mathcal I}}
 \ncm{\U}{{\mathcal U}}
 \ncm{\Hil}{{\mathcal H}}

 \ncm{\Z}{{\mathbb{Z}}}
 \ncm{\eps}{\epsilon}
 \ncm{\ran}{{\rangle}}
 \ncm{\lan}{{\langle}}

 \newtheorem{theo}{Theorem}[section]
\newtheorem{cor}[theo]{Corollary}
\newtheorem{lem}[theo]{Lemma}
\newtheorem{prop}[theo]{Proposition}
\newtheorem{remark}[theo]{Remark}
\newtheorem{definition}[theo]{Definition}
\newtheorem{example}[theo]{Example}

\newenvironment{pf}{{\it Proof.}}{QED\vspace{3mm}}

\title{Topologically irreducible representations of the Banach $*$-algebra associated with a dynamical system}
\author{Aki Kishimoto\footnote{E-mail: akiksmt@r3.ucom.ne.jp} and Jun Tomiyama\footnote{E-mail: juntomi@med.email.ne.jp}}
\date{April 2016}
\begin{document}

\maketitle

\begin{abstract}
We describe (infinite-dimensional) irreducible representations of the crossed product C$^*$-algebra associated with a topological dynamical system (based on $\Z$)  and we show that their restrictions to the underling $\ell^1$-Banach $*$-algebra are not algebraically irreducible under mild conditions on the dynamical system. The above description of irreducible representations has two ingredients, ergodic measures on the space and ergodic extensions for the tensor product with type I factors; the latter which may not have been explicitly taken up before will be explored by examples. A new class of ergodic measures is also constructed for irrational rotations on the circle.


\end{abstract}

\section{Introduction}
Let $X$ be a compact metrizable  space and $\sigma$ a homeomorphism of $X$, which forms a {\em classical dynamical system} $\Sigma=(X,\sigma)$. The corresponding {\em C$^*$-dynamical system} is $(C(X),\alpha)$ where $\alpha$ is the automorphism of the continuous functions $C(X)$ on $X$ defined by  $\alpha(b)(x)=b\sigma^{-1}(x),\, x\in X$ for $b\in C(X)$.

We denote by $\ell^1(\Z,C(X))$ the Banach space of $\ell^1$ functions from $\Z$ into $C(X)$, which is a Banach $*$-algebra denoted by $\ell^1(\Sigma)$ when equipped with a product and a $*$-involution as follows:
$$ (fg)(n)=\sum_k f(k)\alpha^k(g(n-k))
$$
and
$$ f^*(n)=\alpha^n(f(-n))^*.
$$
We denote by $\delta_n\in \ell^1(\Sigma)$ for $n\in\Z$ the function $\delta_n(k)=\delta_{n,k}$ on $\Z$ and by $f\in C(X)$ the function $f\delta_0:\Z\to C(X)$. Thus $\delta_n^*=\delta_{-n}$ and $\delta_1f\delta_{-1}=\alpha(f)$ for $f\in C(X)$.
We denote by $C^*(\Sigma)$ the enveloping C$^*$-algebra of $\ell^1(\Sigma)$, also identified with the crossed product C$^*$-algebra of $C(X)$ by $\alpha$. Note that any topologically irreducible representation of $\ell^1(\Sigma)$ (on a Hilbert space) extends to an irreducible representation of $C^*(\Sigma)$ and that  the universal $C^*$-norm on $\ell^1(\Sigma)$, by which $C^*(\Sigma)$ is defined, is determined by these representations.

Each of $C^*(\Sigma)$ and $\ell^1(\Sigma)$ (as a norm-closed algebra generated by $C(X)$ and $\delta_{\pm1}$) enables us to recover $\Sigma$ and so is as good as the other in this sense. Though $\ell^1(\Sigma)$ looks more tamable with its explicit definition, a close examination on $\ell^1(\Sigma)$ reveals complexity or irregularity as an algebra  which $C^*(\Sigma)$ glosses over in exchange of adopting a representation-friendly {\em intangible} norm. A difference between the two objects seems to most manifestly appear in the case of the simplest example with $X$ a singleton, $C(\T)$ and $\ell^1(\Z)$, where $\T=\R/\Z$ is the dual of $\Z$. Then the convolution algebra $\ell^1(\Z)$ is known to have a non-self-adjoint closed ideal while the closed ideals of $C(\T)$ are all self-adjoint. (This fact is directly translated into a general $\Sigma$ if $\sigma$ has a finite orbit.) Another difference may be found on the lines of Kadison's result: If a representation of a C$^*$-algebra is topologically irreducible, then it is automatically algebraically irreducible (\cite{Kad57}; see also \cite{Sak,Ped}). We naturally expect that a topologically irreducible, infinite-dimensional, representation of $\ell^1(\Sigma)$ is not algebraically irreducible. Indeed this is shown for irreducible representations induced from aperiodic orbits in $X$ (\cite{JTtoappear}). There must be other properties which exhibit a stark difference between these two objects,  deserving thorough investigation but beyond the scope of our present research. Thus we are here confined to the problem of irreducible representations. (See \cite{JST12,JT12,JT13,JTtoappear,Tomiyama,Tom92,Tom} for the ideal structures and some irreducible representations).

We will show this algebraic non-irreducibility for all infinite-dimensional irreducible representations of $\ell^1(\Sigma)$ if $\sigma$ preserves a metric on $X$ which induces the right topology and will give a sufficient condition in other cases. For this purpose we first give a procedure for constructing irreducible representations of $C^*(\Sigma)$ in terms of ergodic $\sigma$-quasi-invariant probability measures on $X$ and some {\em ergodic extensions} of the transformation induced by $\sigma$ (Proposition \ref{irreducible}). Then we give the aforementioned result on algebraic non-irreducibility of representations of $\ell^1(\Sigma)$ (Theorem \ref{main}). In Section 4 we elucidate how ergodic extensions may be possible by examples. Specifically given an ergodic transformation $\sigma$ on a probability space $L^\infty(X)$ we ask a question of whether $\sigma$ can be extended to an ergodic transformation on $L^\infty(X)\otimes M_n$ when $n<\infty$. We manage to give a positive answer in the case of Bernoulli shifts (Proposition \ref{Bernoulli}) and irrational rotations on the circle (Proposition \ref{irrational}) by specifying a certain form of unitaries in $L^\infty(X,\mu)\otimes M_n$ for this extension. We also work on unitary equivalence among those ergodic extensions (Propositions \ref{Bernoulli-equiv}, \ref{Bernoulli-equiv1}, \ref{equiv0} and \ref{equiv1}). But we leave the problem unanswered for general ergodic transformations. Finally we construct a new class of ergodic quasi-invariant probability measures on the circle for an irrational rotation, which is neither atomic nor Lesbegue, where the condition of ergodicity seems to require a specific proof (Proposition \ref{sing-cont}).

\section{Irreducible representations}

Let $\pi$ be an irreducible representation of $C^*(\Sigma)$ and let $\mu$ be a probability measure on $X$ such that $\pi|C(X)$ extends to an isomorphism from $L^\infty(X,\mu)$ onto $\pi(C(X))''$. Let $U=\pi(\delta_1)$, a unitary satisfying $\Ad\,U\pi=\pi\alpha$ on $C(X)$, which implies that $\mu$ must be $\sigma$-quasi-invariant. Since $\pi(C(X))''\cap U'\subset \pi(C(X))'\cap U'=\pi(C^*(\Sigma))'=\C1$, we conclude that $\Ad\,U$ acts on $\pi(C(X))''$ ergodically; thus $\mu$ is ergodic.

\begin{lem}
Let $\pi$ be an irreducible representation of $C^*(\Sigma)$ on a Hilbert space $\Hil_\pi$. Then there is an ergodic $\sigma$-quasi-invariant probability measure $\mu$ on $X$ and a Hilbert space $\Hil$ such that $\Hil_\pi$ is identified with $L^2(X,\mu)\otimes \Hil$ and $\pi(f)=M_f\otimes 1$ for $f\in C(X)$, where $M_f$ denote the multiplication of $f$ on $L^2(X,\mu)$.
\end{lem}
\begin{pf}
Note that $\Hil_\pi$ is separable. The commutant $\pi(C(X))'$ is isomorphic to
$$
\bigoplus_n L^\infty(X,\mu)E_n\otimes \B(\Hil_n)
$$
where $E_n$ is a projection in $L^\infty(X,\mu)$ and $\Hil_n$ is an $n$-dimensional Hilbert space with $n$ including infinity. Since $\Ad\,U$ acts on $\pi(C(X))'$ ergodically, we conclude that $E_n$ must be zero or $1$ and there is only one direct summand. Hence $\pi(C(X))'$ is unitarily equivalent to $L^\infty(X,\mu)\otimes \B(\Hil)$ for some Hilbert space $\Hil$ where $\pi(f)$ corresponds to $M_f\otimes 1$.
\end{pf}

Define a unitary $V$ on $L^2(X,\mu)$ by $(V\xi)(x)=(d\mu\sigma^{-1}/d\mu)^{1/2}(x)\xi\sigma^{-1}(x)$. Then for $f\in C(X)$
$$
(VM_f\xi)(x)=\big(\frac{d\mu\sigma^{-1}}{d\mu}(x)\big)^{1/2}(M_f\xi)(\sigma^{-1}(x))=f(\sigma^{-1}(x))(V\xi)(x),
$$
which implies that $VM_fV^*=M_{\alpha(f)}$.

$\Ad(U(V\otimes1)^*)$ defines an automorphism of $Z=L^\infty(X,\mu)\otimes \B(\Hil)$ which acts trivially on its center. Hence there is a unitary $W\in Z$ such that $\Ad(U(V\otimes1)^*)=\Ad\,W$ on $Z$, i.e., $U(V\otimes1)^*W^*\in Z'\subset Z$ (Proposition 8.9.2 of \cite{Ped}). We may suppose that $U=W(V\otimes1)$ by further modifying $W$ by a central unitary of $Z$ if necessary.

\begin{prop}\label{irreducible}
All the irreducible representations of $C^*(\Sigma)$ are constructed as follows: Choose an ergodic $\sigma$-quasi-invariant probability measure $\mu$ on $X$ and find a unitary $W\in L^\infty(X,\mu)\otimes \B(\Hil)$ for some Hilbert space $\Hil$ such that $\Ad(W(V\otimes 1))$ acts ergodically on $L^\infty(X,\mu)\otimes \B(\Hil)$ where $V$ is the unitary induced by $\sigma$ as above. Then one can define an irreducible representation $\pi$ of $C^*(\Sigma)$ on $L^2(X,\mu)\otimes \Hil$ by
$$
\pi(f)=M_f\otimes 1,\ \ f\in C(X),\ \ \ \ \pi(\delta_1)=W(V\otimes 1).
$$

Let us denote the above representation by $\pi_{(\mu,\Hil,W)}$. Then $\pi_i=\pi_{(\mu_i,\Hil_i,W_i)},\ i=1,2$ are unitarily equivalent with each other if and only if $\mu_1$ and $\mu_2$ are absolutely continuous with each other and $\dim(\Hil_1)=\dim(\Hil_2)$ and there is a unitary operator $\zeta$ from $L^2(X,\mu_2)\otimes \Hil_2$ onto $L^2(X,\mu_1)\otimes \Hil_1$ such that $\pi_1(f)=\zeta\pi_2(f)\zeta^*,\ f\in C(X)$ and $W_1=\zeta W_2\bar{\alpha}(\zeta)^*$, where $\bar{\alpha}(\zeta)=(V_2\otimes1)\zeta (V_1\otimes 1)^*$ and $V_i$ is the $V$ defined for $\mu=\mu_i$.
\end{prop}
\begin{pf}
The first half is proved before this proposition. The unitary equivalence is by definition the existence of $\zeta$ above. The other conditions are redundant but follow from this.
\end{pf}

\begin{prop}\label{noeigenvectors}
If $\pi=\pi_{(\mu,\Hil,W)}$ with $\dim(\Hil)>1$ then $U=\pi(\delta_1)$ has no eigenvectors. Moreover $U$ does not satisfy the equality $U\xi=Y\xi$ for any unit vector $\xi\in L^2(X,\mu)\otimes\Hil$ and any unitary $Y\in L^\infty(X,\mu)\otimes 1$.
\end{prop}
\begin{pf}
Suppose that $U=\pi(\delta_1)$ has a {\em generalized eigenvector}, say $U\xi=Y\xi$ for some unit vector $\xi\in L^2(X,\mu)\otimes\Hil$ and some unitary $Y\in L^\infty(X,\mu)\otimes 1$.  Since
$$(d\mu\sigma^{-1}/d\mu)^{1/2}(x)\xi(\sigma^{-1}(x))= W(x)^*Y(x)\xi(x)\ {\rm a.e.}
$$
we deduce that the set of $x$ with $\|\xi(x)\|=0$ is $\sigma$-invariant. Hence $\xi(x)\not=0$ a.e. Let $e_1(x)=\xi(x)/\|\xi(x)\|,\ x\in X$, which forms a vector $e_1$ in $L^2(X,\mu)\otimes \Hil$. There is a family $e_i=e_i(x),\ i=2,3,\ldots$ of vectors in $L^2(X,\mu)\otimes \Hil$ such that $(e_i(x))_{i\geq 1}$ is a complete orthnormal system in $\Hil$ for almost all $x$ (3.3 of \cite{Sak}). Hence the projection onto the closed subspace $\pi(C(X))\xi$ is a proper projection in the commutant of $\pi(C^*(\Sigma))$, which contradicts the irreducibility of $\pi$. Thus $U$ does not have a {\em generalized eigenvector}.
\end{pf}

If $\mu$ is $\sigma$-invariant and $\dim(\Hil)=1$ in $\pi=\pi_{(\mu,\Hil,W)}$ then $U=\pi(\delta_1)=WV$ satisfies $U1=W1$ where $1$ is regarded as a function in $L^2(X,\mu)$.

\begin{lem}
Suppose that $\mu$ is an ergodic $\sigma$-quasi-invariant probability measure on $X$. Then $\mu$ is either atomic in which case there is $x\in X$ such that $\{\sigma^n(x)\ |\ x\in\Z\}$ has measure 1, or completely non-atomic.
\end{lem}
\begin{pf}
Note that $\mu$ is the sum of an atomic part and a completely non-atomic part and that the decomposition into these parts is unique. Since $\mu$ is ergodic, one of them must be zero. If it is atomic then $\mu$ must be supported by an orbit as it is ergodic.
\end{pf}

\begin{prop}\label{cocyclevanishing}
Let $\pi=\pi_{(\mu,\Hil,W)}$ be an irreducible representation of $C^*(\Sigma)$ and suppose that $\mu$ is atomic. Then $\Hil=\C$. Moreover if $\mu$ has infinite support then $\pi$ is unitarily equivalent to $\pi_{(\mu,\C,1)}$ and if $\mu$ consists of $k$ atoms then $\pi$ is unitarily equivalent to $\pi_{(\mu,\C,\lambda)}$ where $\lambda\in \{e^{2\pi i\theta}\in \C\ |\ 0\leq \theta <1/k\}$.
\end{prop}
\begin{pf}
Suppose that $L^2(X,\mu)\cong \ell^2(\Z)$ and $V$ is the unitary induced by the shift $\sigma: n\mapsto n+1$. We identify $W$ with the sequence $(W_n)_{n\in\Z}$ where $W_n$ is a unitary on $\Hil$. Define a sequence $(\zeta_n)$ of unitaries on $\Hil$ as follows: $\zeta_0=1, \zeta_n=\zeta_{n-1}W_n^*$ for $n>0$, and $\zeta_n=\zeta_{n+1} W_{n+1}$ for $n<0$, and set $\zeta=(\zeta_n)\in \ell^\infty(\Z)\otimes \B(\Hil)$. Then $\zeta_n W_n \zeta_{n-1}^*=1$ for all $n$, which implies that $\zeta W \bar{\alpha}(\zeta)^*=1$, where $\bar{\alpha}=\Ad(V\otimes1)$ is the shift on $\ell^\infty(\Z)\otimes \B(\Hil)$. Thus the fixed point algebra of $\ell^\infty(\Z)\otimes \B(\Hil)$ under $\Ad(W(V\otimes 1))$ is $\Ad\zeta(1\otimes \B(\Hil))$ as $\zeta W(V\otimes 1)\zeta^*=V\otimes 1$. Hence it follows that $\Hil\cong\C$ and $\pi$ is unitarily equivalent to $\pi_{(\mu,\C,1)}$.

Suppose that $L^2(X,\mu)\cong \ell^2(\Z/k\Z)$ with $\Z/k\Z$ identified with $\{0,1,\ldots,k-1\}$ and $V$ is the unitary induced by the shift. Let $Z$ be a unitary in $\B(\Hil)$ such that $Z^k=W_kW_{k-1}\cdots W_1$ where $W=(W_n)_n$ with $W_k=W_0$. Let $\zeta_0=1$ and $\zeta_n=Z\zeta_{n-1}W_n^*$ for $n=1,2,\ldots,k-1$. Then $\zeta_n W_n \zeta_{n-1}^*=Z$ (e.g., $\zeta_0W_0\zeta_{k-1}^*=W_0W_{k-1}\cdots W_1 Z^{-k+1}=Z$ for $n=0$). With $\zeta=(\zeta_n)$ the fixed point algebra of $\ell^2(\Z/k\Z)\otimes \B(\Hil)$ under $\Ad(W(V\otimes1))$ is $\Ad\zeta(1\otimes \B(\Hil)\cap Z')$. Hence it follows that $\Hil\cong\C$. We may assume that $Z$ is a constant as in the statement.
\end{pf}

\begin{remark}
The irreducible representations presented in the above proposition have been already explored in \cite{Tomiyama}. In particular, in the latter case, $\pi_{(\mu,\C,e^{i2\pi\theta})}$, which is equivalent to $\pi_{(\mu,\C, W_\theta)}$ with $W_\theta=(e^{i2\pi k\theta},1,\ldots,1)$, are mutually disjoint for $0\leq \theta <1/k$.
\end{remark}

\section{Topological versus algebraic}

Let $\pi$ be an irreducible representation of $C^*(\Sigma)$. We assume that $\pi=\pi_{(\mu,\Hil,W)}$ as in Proposition \ref{irreducible}. We will show that $\pi|\ell^1(\Sigma)$ is not algebraically irreducible if $L^2(X,\mu)$ is infinite-dimensional under some condition on the quasi-invariance of $\mu$.

\begin{lem}\label{closed-graph}
Let $\pi=\pi_{(\mu,\Hil,W)}$ be as above. Define a bounded linear map $T_\Phi$ of $\ell^1(\Sigma)$ into $L^2(X,\mu)\otimes \Hil$ for a unit vector $\Phi\in L^2(X,\Sigma)\otimes \Hil$ by
$$
T_\Phi(b)=\pi_\mu(b)\Phi.
$$
Let $L_\Phi$ be the kernel of $T_\Phi$, which is a closed left ideal of $\ell^1(\Sigma)$. If $T_\Phi$ is surjective, then there is a constant $K_\Phi>0$ such that
$$
\|b+L_\Phi\|_1\leq K_\Phi\|T(b)\|.
$$
\end{lem}
\begin{pf}
This follows from the closed graph theorem.
\end{pf}

Define $F_k\in L^2(X,\mu)$ for $k\in\Z$ by
$$
F_k(x)=(\frac{d\mu\sigma^{-k}}{d\mu}(x))^{1/2}
$$
and note that $V$ satisfies that $(V^k\xi)(x)=F_k(x)\xi(\sigma^{-k}(x)),\ \xi\in L^2(X,\mu)$.

\begin{lem}\label{norm-estimate}
Let $\eta$ be a unit vector of $\Hil$. Let $A$ be a measurable subset of $X$ with $\mu(A)>0$ and let $S=\sum_k a_k\delta_k\in \ell^1(\Sigma)$ be such that
$$
\pi_\mu(S)1\otimes\eta =\frac{\chi_A}{\mu(A)^{1/2}}\otimes \eta,
$$
where $\chi_A$ is the characteristic function of $A$. Then it follows that $1\leq \sum_k \|a_k\|\mu\sigma^{-k}(A)^{1/2}\leq \|S\|_1\sup_k\mu(\sigma^k(A))^{1/2}$.
\end{lem}
\begin{pf}
Note that $\pi(\delta_k)= (W(V\otimes 1))^k=W_k(V^k\otimes 1)$ where $W_0=1$, $W_k=W\bar{\alpha}(W_{k-1})$ for $k>0$ and $W_k=\bar{\alpha}^{-1}(W^* W_{k+1})$ for $k<0$, and $\bar{\alpha}=\Ad(V\otimes 1)$; $W_k$ are all unitaries in $L^\infty(X,\mu)\otimes\B(\Hil)$. Let $\xi=\chi_A/\mu(A)^{1/2}$, a unit vector. Then we compute:
$$
1=\sum_k\lan \xi\otimes \eta, (a_k\otimes 1)W_k(F_k\otimes \eta)\ran=\sum_k \lan \xi\otimes\eta, (a_k\otimes 1)W_k(\chi_A F_k)\otimes \eta \ran
$$
which is at most $\sum_k\|a_k\|\| \chi_A F_k\otimes \eta\|=\sum_k\|a_k\|\mu\sigma^{-k}(A)^{1/2}$.
\end{pf}

\begin{lem}
Suppose that $\mu$ is non-atomic and that there is a metric $d$ on $X$ such that $d$ induces the topology on $X$ and $\sigma$ preserves $d$, i.e., $d(x,y)=d(\sigma(x),\sigma(y))$ for $x,y\in X$. Then for any $\epsilon>0$ there is a measurable subset $A$ of $X$ such that $0<\sup_k\mu(\sigma^k(A))<\epsilon$.
\end{lem}
\begin{pf}
Suppose, on the contrary, that there is an $\epsilon>0$ such that $\sup_k\mu(\sigma^k(A))>\epsilon$ for any $A$ with $\mu(A)>0$. Let $x\in X$ be such that any open neighborhood of $x$ has positive measure. Let $U_n=\{y\in X\ |\ d(x,y)<1/n\}$ for $n\in\N$. Then there is a $k_n\in \Z$ such that $\mu(\sigma^{k_n}(U_n))>\epsilon$. Since $X$ is compact there is a subsequence in $(\sigma^{k_n}(x))$ converging, say to $z\in X$. Then it follows that any open neighborhood of $z$ contains $\sigma^{k_n}(U_n)$ for some $n$ and hence has measure greater than $\epsilon$, which implies $\mu(\{z\})\geq \epsilon$. Hence $\mu$ is atomic, a contradiction.
\end{pf}

Let $\Lambda=\{\sup_k\mu(\sigma^{k}(A))\ |\ A\ {\rm measurable\ with}\ \mu(A)>0  \}$, which is a subset of $(0,1]$. If $\mu$ is $\sigma$-invariant and is non-atomic, then $\Lambda=(0,1]$. If $\mu$ is atomic, then $\Lambda\subset [\lambda_0,1]$ where $\lambda_0=\sup_k\mu\sigma^k(A)>0$ with $A$ an atom.

\begin{theo}\label{main}
Let $\pi=\pi_{(\mu,\Hil,W)}$ be as above and suppose that $\mu$ is non-atomic. Assume that for any $\epsilon>0$ there is a measurable subset $A$ of $X$ such that $0<\sup_k\mu\sigma^k(A)<\epsilon$ $($which follows if $\mu$ is $\sigma$-invariant or $\sigma$ preserves a metric on $X$ which induces the topology$)$. Then $\pi|\ell^1(\Sigma)$ is not algebraically irreducible.
\end{theo}
\begin{pf}
The parenthesized statement follows from the previous lemma and the remark before this theorem.

Suppose that $\pi|\ell^1(\Sigma)$ is algebraically irreducible. Let $\eta$ be a unit vector of $\Hil$. Then the map $T_{1\otimes \eta}$ from $\ell^1(\Sigma)$ into $L^2(X,\mu)\otimes \Hil$ is surjective. Hence Lemma \ref{closed-graph} gives a constant $K>0$ satisfying: For any unit vector $\Psi\in L^2(X,\mu)$ there is an $S\in\ell^1(\Sigma)$ such that $\pi(S)1\otimes \eta=\Psi$ and $\|S\|_1\leq K$. Lemma \ref{norm-estimate} leads us to a contradiction under the hypothesis by taking $\Psi=\chi_A/\mu(A)^{1/2}\otimes \eta$ for $A$ with small $\sup_k\mu\sigma^k(A)$.
\end{pf}

\begin{cor}
Let $\Sigma=(X,\sigma)$ and suppose that $X$ is a metric space and $\sigma$ preserves the metric on $X$ and has no periodic points. Then $\pi|\ell^1(\Sigma)$ is not algebraically irreducible for any irreducible representation $\pi$ of $C^*(\Sigma)$.
\end{cor}
\begin{pf}
Under the hypothesis all irreducible representations are infinite-dimensional.

If $\mu$ is completely non-atomic then this follows from Theorem \ref{main}. If $\mu$ is atomic, then $L^2(X,\mu)\cong \ell^2(\Z)$ and this is proved in \cite{JTtoappear}.

Let us repeat the proof in the atomic case, which seems subtler, but simpler, than the one of Theorem \ref{main}. In this case we may work in $\ell^2(\Z)$ with $\sigma$ the shift on $\Z$. Denote by $\xi_n$ the function in $\ell^2(\Z)$ defined by $\xi_n(k)=\delta_{n,k}$. Suppose that there is $S=\sum_k a_k\delta_k\in\ell^1(\Sigma)$ such that $\pi(S)\xi_0=\sum_{k=1}^\infty k^{-1}\xi_k$. Since $\pi(S)\xi_0=\sum_k \pi(a_k)\xi_k$ it follows that $\pi(a_k)\xi_k=k^{-1}\xi_k$ for $k=1,2,\ldots$, which implies that $\|a_k\|\geq 1/k$ for $k\geq 1$. This contradicts $\sum_k\|a_k\|<\infty$.
\end{pf}


\section{Ergodic extensions}

The observation on irreducible representations of $C^*(\Sigma)$ in Proposition \ref{irreducible} gives rise to a problem of whether given an ergodic transformation $\sigma$ on $(X,\mu)$ there is a unitary $W\in Z=L^\infty(X,\mu)\otimes \B(\Hil)$ such that $\gamma=\Ad(W(V\otimes 1))$ acts on $Z$ ergodically for a given Hilbert space $\Hil$. Here $V$ is the unitary on $L^2(X,\mu)$ defined by $(V\xi)(x)=(d\mu\sigma^{-1}/d\mu)^{1/2}(x)\xi(\sigma^{-1}(x))$, which implements $\alpha$ on $L^\infty(X,\mu)$ (where $\alpha(f)=f\sigma^{-1},\ f\in L^\infty(X,\mu)$ as before). We have shown that if $\mu$ is atomic and $\dim (\Hil)>1$ then there is no such $W$. Hence we assume that $\mu$ is non-atomic and call this the problem of {\em ergodic extensions}. We shall write $V$ in place of $V\otimes 1$ from now on.

Let $W=\int_X^\oplus W(x)d\mu(x)\in Z$ where we assume that $W(x)$ is a unitary on $\Hil$ for all $x\in X$ (as $W$ is a unitary). Let $T=\int_X^\oplus T(x)d\mu(x)\in Z$. Then $T\in Z^\gamma=Z\cap (WV)'$ if and only if
$$
T\sigma^{-1}(x)=W(x)^*T(x)W(x)
$$
almost everywhere. In particular $\|T\sigma^{-1}(x)\|=\|T(x)\|$ a.e., which implies that $\|T(x)\|=\|T\|$ almost everywhere.  If $\Hil$ is finite-dimensional it then easily follows:

\begin{lem}
When $\dim \Hil=n<\infty$, $Z\cap (WV)'$ is isomorphic to a $*$-subalgebra of $\B(\Hil)$ for any $W$.
\end{lem}
\begin{pf}
Let $T\in Z\cap (WV)'$. Then it follows that $\det(T(x)-\lambda1)$, as a polynomial of order $n$ in $\lambda$, is almost constant. Hence $T$ has at most $n$ eigenvalues, which implies that $Z\cap (WV)'$ is finite-dimensional. Hence one can show that there is an $x\in X$ such that $Z\cap (WV)'\ni T\mapsto T(x)$ is an injective homomorphism.
\end{pf}

The problem we cannot answer in general is whether there is a unitary $W\in Z=L^\infty(X,\mu)\otimes\B(\Hil)$ for $\Hil\not=\C1$ such that $Z\cap (WV)'=\C1$ when $\mu$ is non-atomic, let alone how to classify those $W$ modulo unitary equivalence.

In what follows we restrict ourselves to two specific examples of $\Sigma=(X,\sigma)$; Bernoulli shifts and irrational rotations. If $\Sigma$ is a Bernoulli shift then $C^*(\Sigma)$ is not simple but primitive and if $\Sigma$ is an irrational rotation on $X=\T$ then $C^*(\Sigma)$ is simple (see \cite{Ped,Tomiyama}). Among many ergodic measures on $X$ we shall choose specific invariant probability measures on $X$ and discuss ergodic extensions.

\medskip

Let $\Lambda$ be a finite set of more than one elements. Let $X=\Lambda^\Z=\{(x_n)\ |\ x_n\in \Lambda\}$ and $\sigma$ be the shift on $X$ to the right and define an automorphism $\alpha$ on $C(X)$ by $\alpha(f)(x)=f\sigma^{-1}(x)$. Then the fixed point algebra $C(X)^\alpha$ of $C(X)$ under $\alpha$ is $\C1$. This follows because there is an $x\in X$ whose $\sigma$-orbit is dense in $X$.

First we consider a $C(X)$-version instead of $L^\infty(X,\mu)$, which is considerably simpler, i.e., we assert that for any integer $n>1$ there is an automorphism $\beta$ of $C(X)\otimes M_n$ such that $\beta(f\otimes 1)=\alpha(f)\otimes 1$ for $f\in C(X)$ and $(C(X)\otimes M_n)^\beta=\C1$.

We shall prove this assertion. Define a diagonal unitary $u$ by $u=1\oplus \omega\oplus \omega^2\oplus \cdots \oplus \omega^{n-1}$ with $\omega=e^{2\pi i/n}$ and a shift unitary $v\in M_n$ such that
$$
vuv^*=\omega u.
$$
Note that $M_n\cap \{u,v\}'=\C1$. Let $C_1$ be a non-empty proper subset of $\Lambda$ and let $C=\{x\in X\ |\ x_0\in C_1\}$, a closed and open subset of $X$. Let $D=X\setminus C=\{x\in X\ |\ x_0\in \Lambda\setminus C_1\}$. Define a unitary $W\in C(X)\otimes M_n$ by
$$
W=\chi_C\otimes u^*+\chi_D\otimes v^*
$$
and an automorphism $\beta$ of $C(X)\otimes M_n$ by $\beta=\Ad\,W\circ (\alpha\otimes\id)$. Note that for $f\in C(X)\otimes M_n$,
$$
\beta(f)(x)=W(x)f(\sigma^{-1}(x))W(x)^*=\begin{cases} \Ad\,u^*(f(\sigma^{-1}(x)))& x\in C\\
                                                      \Ad\,v^*(f(\sigma^{-1}(x)))& x\in D.
                                                               \end{cases}
$$
Hence if $\beta(f)=f$ then
$$
f\sigma^{-1}(x)=\Ad\,W(x)^*(f(x)).
$$

\begin{lem}\label{values}
Let $T\in (C(X)\otimes M_n)^\beta$. Then there is a $T_0\in M_n$ such that $T$ takes all the values in $\Gamma=\{\Ad(u^pv^q)(T_0); p,q=0,1,2,\ldots,n-1\}$.
\end{lem}
\begin{pf}
Since $u^n=1=v^n$ and $\Ad\,u\Ad\,v=\Ad\,v\Ad\,u$, the finite subset $\Gamma$ of $M_n$ is invariant under $\Ad\,u$ and $\Ad\,v$ for any $T_0\in M_n$.

Let $T\in Z\cap (WV)'$ and $x^0\in X$ be such that the orbit $\{\sigma^k(x^0)\ |\ k\in\Z\}$ is dense in $X$. Let $T_0=T(x^0)$. Then $T\sigma^{-1}(x^0)=\Ad\,u(T_0)$ or $\Ad\,v(T_0)$ depending on $x^0\in C$ or $x^0\in D$; thus $T\sigma^{-1}(x^0)\in \Gamma$. Repeating this process it follows that $T\sigma^k(x^0)\in \Gamma$ for all $k$. Since $T$ is continuous on $X$ we deduce that $T(x)\in \Gamma$ for all $x\in X$.

Let $c(x,k)=\#\{i\ |\ x_i\in C_1, 0\leq i<k\}$ for $x\in X$ and $k\in\N$ and $d(x,k)=k-c(x,k)$, where $\# B$ denote the number of points in a set $B$. Then we obtain that $T(\sigma^{-k}(x))=\Ad(u^{c(x,k)}v^{d(x,k)})(T(x))$ for $k\in\N$.

Let $S$ be the cylinder subset of $X$ consisting of $x$ with $x_0,x_1,\ldots,x_{n^2-1}$ specified as follows: $x_i\in C_1$ for $0\leq i<n$, $x_n\not\in C_1$, $x_{n+1+i}\in C_1$ for $0\leq i<n-1$, $x_{2n}\not\in C_1$, $x_{2n+1+i}\in C_1$ for $0\leq i<n-1$, $x_{3n}\not\in C_1, \ldots, x_{(n-1)n}\not\in C_1, x_{(n-1)n+1+i}\in C_1$ for $0\leq i<n-1$. Then for $x\in S$ the set of pairs $(c(x,k)+n\Z,d(x,k)+n\Z)$ with $k=0,1,2,\ldots,n^2-1$ exhausts the whole $\Z/n\Z\times \Z/n\Z$. This shows that $T(\sigma^k(x)), k=0,1,2,\ldots,n^2-1$ exhausts the whole $\Gamma$.
\end{pf}

By the above lemma $T\in (C(X)\otimes M_n)^\beta$ takes a finite number of values, say $T_i,\ i=1,2, \ldots,m$, where $m$ divides $n^2$. If $m=1$ this implies that $\Ad(u^pv^q)(T_0)=T_0$ for all $p,q$, i.e., $T_0\in\C1$ or $T\in 1\otimes \C1$. So we assume that $m>1$.

Note that $F_i=T^{-1}(T_i)$ is a closed and open subset of $X$. Then $F_i$ is a cylinder subset. (Given $z\in F_i$ there is a (open) cylinder subset $U(z)$ such that $z\in U(z)\subset F_i$. Since $F_i$ is compact one can find a finite number of $U(z)$ whose union  equals $F_i$.) Hence there is an $N\in \N$ such that all $F_i$'s are determined by subsets of $\prod_{k=-N}^{N-1} \Lambda$. Let $S_i$ be a subset of $\prod_{k=-N}^{N-1}\Lambda$ such that $F_i=\{x\in X\ |\ (x_{-N},x_{-N+1},\ldots, x_{N-1})\in S_i\}$.

Let $y\in S_1$ and $z\in S_2$. We shall construct an element $x\in X$ containing $y,z$ as segments whose existence gives a contradiction.

Let $m_1=\#\{i\ |\ y_i\in C_1, 0\leq i\leq N-1\}$ and $m_2=\#\{i\ |\ z_i\in C_1, -N\leq i\leq -1\}$ and let $\ell_1=N-m_1$ and $\ell_2=N-m_2$. Let $a,b$ be integers between $0$ and $n-1$ such that $m_1+m_2+a=0, \ell_1+\ell_2+b=0$ modulo $n$.  We define $x\in X$ as an element satisfying the following conditions:
\begin{align*}
x_i&=y_i,  &-N\leq i\leq N-1,\\
x_{i}&\in C_1, & N\leq i\leq N+a-1,\\
x_{i}&\not\in C_1, & N+a\leq i\leq N+a+b-1, \\
x_{i}&=z_{i-2N-a-b},& N+a+b\leq i\leq 3N+a+b-1.
\end{align*}
Then $x\in F_1$ and $\sigma^{-(2N+a+b)}(x)\in F_2$ (as $\sigma^{-(2N+a+b)}(x)_i=x_{2N+a+b+i}=z_{i}$ for $-N\leq i\leq N-1$). Since $c(x,2N+a+b)=m_1+a+m_2=0$ (mod $n$) and $d(x,2N+a+b)=\ell_1+b+\ell_2=0$ (mod $n$), it follows that
$$
T(\sigma^{-(2N+a+b)}x)=\Ad(u^{c(x,2N+a+b)}v^{d(x,2N+a+b)})(T_1)=T_1,
$$
which contradicts that $T(\sigma^{-(2N+a+b)}(x))=T_2$ following from $\sigma^{-(2N+a+b)}(x)\in F_2$. Thus one can conclude that $m=1$.

\begin{prop}
Let $X=\Lambda^\Z$ and $\sigma$ the shift on $X$ as above. If $\alpha$ is the automorphism of $C(X)$ induced by $\sigma$ then $C(X)^\alpha=\C1$. If $n$ is an integer greater than 1 and $\beta=\Ad\,W(\alpha\otimes 1)$ is an automorphism of $C(X)\otimes M_n$ with $W$ as above, it follows that $(C(X)\otimes M_n)^\beta=\C1$.
\end{prop}

We will now prove the $L^\infty$-version of the above result. Let $X=\Lambda^\Z$ and $\sigma$ the shift on $X$ as above. Let $\mu_1$ be a probability measure on $\Lambda$ such that $\mu_1(\{\lambda\})>0$ for all $\lambda\in \Lambda$ and define a measure $\mu$ on $X$ as the infinite product of copies of $\mu_1$. Then $\mu$ is a $\sigma$-invariant probability measure on $X$. Define a unitary $V$ on $L^2(X,\mu)$ as the unitary induced by $\sigma$ as before. Then $\Ad\,V$ acts on $L^\infty(X,\mu)$ ergodically.

This is shown in a standard way. For any pair $A,B$ of cylinder subsets of $X$ we obtain that $\mu(A\cap \sigma^k(B))\to \mu(A)\mu(B)$ as $k\to\infty$. It then follows that this is true for any measurable subsets $A,B$. If $A$ is a $\sigma$-invariant subset, i.e., $\mu(A\setminus \sigma(A))=0=\mu(\sigma(A)\setminus A)$ then it follows that $\mu(A)=\mu(A\cap \sigma^k(A))=\mu(A)^2$, i.e., $\mu(A)=0$ or $1$. Hence $\sigma$ is ergodic.

Next we will prove: For any $n>1$ there is a unitary $W$ in $Z=L^\infty(X,\mu)\otimes M_n$ such that $\Ad(WV)$ acts on $Z$ ergodically.

We have defined the unitaries $u,v\in M_n$ and define a unitary $W\in L^\infty(X,\mu)\otimes M_n$ as before:
$$
W=\chi_C\otimes u^*+\chi_D\otimes v^*,
$$
where $C=\{x\in X\ |\ x_0\in C_1\}$ and $D=X\setminus C$.

\begin{lem}\label{finite}
Let $T\in L^\infty(X,\mu)\otimes M_n\cap (WV)'$. Then there is a $T_0\in M_n$ such that $T$ takes values in $\Gamma=\{\Ad(u^pv^q)(T_0)\ |\  p,q=0,1,2,\ldots,n-1\}$. Moreover $T$ takes all the values in $\Gamma$.
\end{lem}
\begin{pf}
As we have remarked, the subset $\Gamma$ of $M_n$ is invariant under $\Ad\,u$ and $\Ad\,v$ for any $T_0\in M_n$.

Let $T\in Z\cap (WV)'$. Let $T_0\in M_n$ be such that $\|T_0\|=\|T\|$ and $\{x\in X\ |\ \|T(x)-T_0\|<\epsilon \}$ has positive measure for any $\epsilon>0$. Let $B_\epsilon=\{S\in M_n\ |\ \|S-T_0\|<\epsilon\}$. Then it follows that $T(x)\in F_\epsilon=\bigcup_{p,q}\Ad(u^pv^q)(B_\epsilon)$ for almost all $x$. (If $T(x)\in F_\epsilon$ then $T\sigma^{-1}(x)\in F_\epsilon$ since $T\sigma^{-1}(x)=\Ad\,u(T(x))$ or $\Ad\,v(T(x))$. Repeating this it follows that $T\sigma^k(x)\in F_\epsilon$ for all $k$.) By taking the intersection of $F_\epsilon$ with $\epsilon>0$ we conclude that $\{x\in X \ |\ T(x)\in \Gamma\}$ has full measure.

Let $X_{p,q}=\{x\in X\ |\ T(x)=\Ad(u^pv^q)(T_0)\}$ for $p,q=0,1,\ldots,n-1$. Let $S$ be the cylinder subset defined in the proof of Lemma \ref{values}. There is a pair $p,q$ such that $S\cap X_{p,q}$ is not a null set. Then it follows from the property of $S$ that $\sigma^{-k}(S\cap X_{p,q}),\ k=0,1,2,\ldots,n^2-1$ visit all $X_{i,j}$. (For example $\sigma^{-k}(S\cap X_{p,q})\subset X_{p+k,q}$ for $k=0,1,\ldots,n-1, \sigma^{-n}(S\cap X_{p,q})\subset X_{p+n-1,q+1}$, etc.) This proves the last statement.
\end{pf}

If $\Gamma$ is a singleton, then $T_0\in \C1$ because $u,v$ generate the whole $M_n$. Suppose that $\Gamma$ includes at least two points; say $T_0$ and $T_1=\Ad(u^{p'}v^{q'})(T_0)\not=T_0$ for some $p',q'$. Let $F_i=\{x\in X\ |\ T(x)=T_i\}$ for $i=0,1$ and let $\eta=\min\{\mu(F_i)\ |\ i=0,1\}>0$. Let $K_i$ be a compact subset such that $K_i\subset F_i$ and $\mu(K_i)>\eta/2$. Let $\epsilon>0$ be very small and choose a cylinder subset $O_i$ such that $K_i\subset O_i$ and $\mu(O_i\setminus K_i)<\epsilon$. (We will specify $\epsilon>0$ later.) We choose $N\in\N$ such that $O_i$ can be regarded as a subset of $\prod_{-N}^{N-1}\Lambda$ for $i=0,1$. We assume that $N$ is a multiple of $n$.

Let
\begin{align*}
G_{p}=& \{x\in X\ |\ \sum_{i=0}^{N-1}\chi_{C_1}(x_i)=p\ {\rm mod}\ n\}, \\
H_p=& \{x\in X\ |\ \sum_{i=-N}^{-1}\chi_{C_1}(x_i)=p\ {\rm mod}\ n\},\\
L_p=&\{x\in X\ |\ \sum_{i=N}^{N+n-1} \chi_{C_1}(x_i)=p\ {\rm mod} \ n\}
\end{align*}
for $p=0,1,\ldots,n-1$. Regarding $O_0\cap G_p$ as a subset of $\prod_{i=-N}^{N-1}\Lambda$ and $L_{n-p-q}$ as a subset of $\prod_{i=N}^{N+n-1}\Lambda$ and $O_1\cap H_p$ as a subset of $\prod_{i=N+n}^{3N+n-1}\Lambda$ (after being translated $2N+n$ to the right), we construct a cylinder set $F$ (corresponding to a subset of  $\prod_{i=-N}^{3N+n-1}\Lambda$) by concatenating triple finite sequences in
$$
\bigcup_{p,q}O_0\cap G_p\times L_{n-p-q}\times O_1\cap H_q.
$$
Then $F\subset O_0$ and $\sigma^{-2N-n}(F)\subset O_1$ (because $O_1\cap H_q$ to be concatenated has been shifted to the right by $2N+n$). Since $\bigcup_{p=0}^{n-1} G_p=X$ etc., we estimate
$$
\mu(F)=\sum_{p}\mu(O_0\cap G_p)\sum_q\mu(L_{n-p-q})\mu(O_1\cap H_q)\geq \eta_0\mu(O_0)\mu(O_1)\geq \eta_0\eta^2/4
$$
where $\eta_0=\min_p \mu(L_p)>0$. We assume that $\epsilon<\eta_0\eta^2/8$ (as $\eta_0$ depends only on $\mu$ and $n$); then it follows that $\mu(F\cap K_0\cap \sigma^{2N+n}(K_1))$ is positive because it is bounded below by
$$
\mu(F)-\mu(F\setminus K_0)-\mu(F\setminus \sigma^{2N+n}(K_1))\geq \mu(F)-2\epsilon>0.
$$

Let $x\in F\cap K_0\cap \sigma^{2N+n}(K_1)$. Then $T(\sigma^{-2N-n}(x))=T(x)$ because $\sum_{i=0}^{2N+n-1}\chi_{C_1}(x_i)=0\ {\rm mod}\ n$ (by the construction of $F$) and $2N+n=0\ {\rm mod}\ n$. This is a contradiction because  $T(\sigma^{-2N-n}(x))=T_1\not=T_0=T(x)$ follows from $\sigma^{-2N-n}(x)\in K_1\subset F_1$ and $x\in K_0\subset F_0$. Thus $\Gamma$ must be a singleton.

\begin{prop}\label{Bernoulli}
Let $X=\Lambda^\Z$ and $\sigma$ the shift on $X$ and $\mu=\prod_\Z\mu_1$ a probability measure on $X$ as above. Then the automorphism $\alpha$ on $L^\infty(X,\mu)$ induced by $\sigma$ satisfies $L^\infty(X,\mu)^\alpha=\C1$. If $n$ is an integer greater than 1 and $\beta=\Ad\,W(\alpha\otimes 1)$ is an automorphism of $L^\infty(X,\mu)\otimes M_n$ with $W$ as above then it follows that $(L^\infty(X,\mu)\otimes M_n)^\beta=\C1$.
\end{prop}

We have defined the subset $C$ by specifying $C_1\subset \Lambda$: $C=\{x\in X\ |\ x_0\in C_1\}$ with $C_1\not=\emptyset,\Lambda$. Let $C_1'$ be anther subset of $\Lambda$ and define $C'=\{x\in X\ |\ x_0\in C_1'\}$ and $D'=X\setminus C'$. Let $W'\in L^\infty(X,\mu)\otimes M_n$ be the corresponding unitary $\chi_{C'}\otimes u^*+\chi_{D'}\otimes v^*$. We consider the problem of when  $WV$ and $W'V$ are unitarily equivalent.

Suppose that there is a unitary $\zeta\in L^\infty(X,\mu)\otimes M_n$ such that $WV=\zeta W'V\zeta^*$, i.e., $W(x)=\zeta(x)W'(x)\zeta(\sigma^{-1}(x))^*$ for almost all $x$ or $\zeta\sigma^{-1}(x)=W(x)^*\zeta(x)W'(x)$.  Then depending on $x_0\in \Lambda$ we have the following cases:
\begin{enumerate}
\item If $x\in C\cap C',\ \ \zeta(\sigma^{-1}(x))=u\zeta(x)u^*$,
\item If $x\in D\cap D',\ \ \zeta(\sigma^{-1}(x))=v\zeta(x)v^*$,
\item If $x\in C\cap D',\ \ \zeta(\sigma^{-1}(x))=u\zeta(x)v^*$,
\item If $x\in D\cap C',\ \ \zeta(\sigma^{-1}(x))=v\zeta(x)u^*$.
\end{enumerate}
Hence if $\zeta(x)$ is defined $\zeta(\sigma^{-1}(x))$ is obtained by applying one of the four maps on $\zeta(x)$:
$$
\phi_1=L_{u}R_{u^*},\ \phi_2=L_{v}R_{v^*},\ \phi_3=L_{u}R_{v^*},\ \phi_4=L_{v}R_{u^*}
$$
depending on $x_0$, where $L_b$ denotes the left multiplication of $b\in M_n$ etc. They satisfy $\phi_i^n=\id$ for $i=1,2,3,4$, $\phi_1\phi_2=\phi_2\phi_1$, and
$$
\begin{array}{lll}
\ &\phi_1\phi_3=\omega\phi_3\phi_1,  & \phi_1\phi_4=\omega^{-1}\phi_4\phi_1,\\
  & \phi_2\phi_3=\omega \phi_3\phi_2,& \phi_2\phi_4=\omega^{-1}\phi_4\phi_2,\\
  &\phi_3\phi_4= \omega^{-2}\phi_4\phi_3,
\end{array}
$$
by the commutation relation $vu=\omega uv$. Note also that $\phi_3^k(\zeta)=\Ad\,u^{k}(\zeta)u^{k}v^{-k}$ and $\phi_4^k(\zeta)=\Ad\,v^{k}(\zeta)v^{k}u^{-k}$.

\begin{lem}
In the above situation there is a unitary $\zeta_0\in M_n$ such that $\zeta$ takes values in $\Gamma=\{\Ad(u^pv^q)(\zeta_0)\omega^k(uv^*)^\ell\ |\ p,q,k,\ell=0,1,2\ldots,n-1\}$ almost everywhere and $\{x\in X\ |\ \zeta(x)=\zeta_0\}$ is not a null set.
\end{lem}
\begin{pf}
Since $\Gamma$ is invariant under $\phi_i,\ i=1,2,3,4$ this can be proved in the same way as Lemma \ref{finite}.
\end{pf}

Suppose that $\zeta$ takes at least two values, $\zeta_0$ above and $\zeta_1\in \Gamma$. Let $F_i=\{x\in X\ |\ \zeta(x)=\zeta_i\}$ for $i=0,1$ and $\eta=\min\{\mu(F_0),\mu(F_1)\}>0$. Let $K_i\subset F_i$ be a compact subset such that $\mu(K_i)>\eta/2$. For any $\epsilon>0$ there are cylinder subsets $O_i$ such that $K_i\subset O_i$ and $\mu(O_i\setminus K_i)<\epsilon$. There is an $N\in\N$ such that each $O_i$ is determined by a subset of $\prod_{i=-N}^{N-1}\Lambda$. We may assume that $N$ is a multiple of $n$.

We have assumed $C\not=C'$: There are six cases depending on which $C\cap C', D\cap D',C\cap D',\D\cap C'$ are empty. ($i$) $C\cap C'=D\cap D'=\emptyset$ ($ii$) only $C\cap C'=\emptyset$, ($iii$) only $C\cap D'=\emptyset$, ($iv$) only $D\cap C'=\emptyset$, ($v$) only $D\cap D'=\emptyset$, ($vi$) none of the above intersections are empty. We will consider each case separately.

Suppose ($i$), i.e., $C=D',D=C'$. In this case only $\phi_3$ and $\phi_4$ appear when we express $\zeta\sigma^{-1}(x)$ in terms of $\zeta(x)$.

We furthermore assume that $n$ is even. Let $x\in O_0$. Then $\zeta\sigma^{-N}(x)$ can be uniquely expressed as $\omega^{2k}\phi_3^\ell\phi_4^m(\zeta(x))$, depending on whether each of $x_0,x_1,\ldots,x_{N-1}$ falls into $C$ or $D$, with $0\leq k\leq n/2-1,\ 1\leq \ell\leq n$, and $n\leq m\leq 2n-1$. Then, since $n+1\leq \ell +m\leq 3n-1$ and $\ell+m=0$ (mod $n$), it follows that $\ell+m=2n$. Hence we can also express this as $\zeta\sigma^{-N}(x)=\phi_3^{\ell-1}\phi_4^k\phi_3\phi_4^{m-k}(\zeta(x))$. We denote by $O_0(k,\ell,m)$ the set of $x\in O_0$ with $x_0,x_1,\ldots,x_{N-1}$ giving this expression on $\zeta\sigma^{-N}(x),\zeta(x)$. The union of $O_0(k,\ell,m)$ with all possible $k,\ell,m$ equals $O_0$.

Let $y\in O_1$. Then in the same way $\zeta(y)=\phi_3^{\ell'-1}\phi_4^{k'}\phi_3\phi_4^{m'-k'}(\zeta\sigma^{N}(y))$ depending on $y_{-N},y_{-N+1},\ldots,y_{-1}$, for $k',\ell',m'$ with $0\leq k'\leq n/2-1,\ 1\leq \ell'\leq n$, and $n\leq m'\leq 2n+1$. We denote by $O_1(k',\ell',m')$ the set of $y\in O_1$ with $y_{-N},y_{-N+1},\ldots,y_{-1}$ giving rise to this relation on $\zeta(y), \zeta\sigma^N(y)$.

Let $L(k,\ell,m)=C^{m-1}\times D^k\times C\times D^{\ell-k}\subset \Lambda^{2n}$. We define $F$ to be the cylinder set of $X$ corresponding to
$$
\bigcup_{k,\ell,m,k',\ell',m'}O_0(k,l,m)\times L(k,\ell,m)\times L(k',\ell',m')\times O_2(k',\ell',m')
$$
as a subset of $\prod_{i=-N}^{3N+4n-1}\Lambda$. Then we deduce for $x\in F$ that $\zeta(\sigma^{-2N-4n}(x))=\zeta(x)$ since
$$
\zeta(\sigma^{-N-2n}(x))=\phi_4^{\ell-k}\phi_3\phi_4^k\phi_3^{m-1}\phi_3^{\ell-1}\phi_4^k\phi_3\phi_4^{m-k}(\zeta(x))
=\zeta(x)
$$
(because $\phi_3^{m-1}\phi_3^{\ell-1}=\phi_3^{-2}$, $\phi_3\phi_4^k=\omega^{-2k}\phi_4^k\phi_3$, $\phi_4^k\phi_3=\omega^{2k}\phi_3\phi_4^k$ etc.)
and
$$\zeta(\sigma^{-2N-4n}(x)) =\phi_3^{\ell'-1}\phi_4^{k'}\phi_3\phi_4^{m'-k'}\phi_4^{\ell'-k'}\phi_3\phi_4^{k'}\phi_3^{m'-1}(\zeta(\sigma^{-N-2n}(x)))=\zeta(\sigma^{-N-2n}(x)).
$$
Note that $\mu(F)\geq \mu(O_0)\mu(O_1)\eta_0^2>\eta_0\eta^2/4$ where $\eta_0=\min\{\mu(L(k,\ell,m))\}>0$ is independent of the choice of $K_i$ etc. If $\epsilon<(\eta_0\eta)^2/8$ then $K_0\cap \sigma^{2N+4n}(K_1)\cap F$ is not a null set, which contradicts $\zeta(\sigma^{-2N-4n}(x))=\zeta(x)$ on $F$ as shown as before. Hence $\zeta(x)$ takes just one value $\zeta_0$ almost everywhere. Thus we conclude that $\zeta_0=u\zeta_0v^*$ and $\zeta_0=v\zeta_0 u^*$ or
$$
u=\zeta_0 v \zeta_0^*,\ \ v=\zeta_0u\zeta_0^*.
$$
So the pair $(u,v)$ maps to $(v,u)$ under $\Ad\,\zeta_0$; this happens when and only when $n=2$. In the case $n=2$ we may take
$$
\zeta_0=\frac{1}{\sqrt{2}}\begin{pmatrix} 1& 1\\ 1& -1 \end{pmatrix}.
$$

We now assume that $n$ is odd. Given $x\in O_0$ the value $\zeta\sigma^{-N}(x)$ can be uniquely expressed as $\omega^{2k}\phi_3^\ell\phi_4^m(\zeta(x))=\phi_3^{\ell-1}\phi_4^k\phi_3\phi_4^{m-k}(\zeta(x))$ with $0\leq k\leq n-1,\ 1\leq \ell\leq n$, and $n\leq m\leq 2n-1$. (Then $\ell+m=2n$.) We denote by $O_0(k,\ell,m)$ the set of $x\in O_0$ with $x_0,x_1,\ldots,x_{N-1}$ giving this expression on $\zeta\sigma^{-N}(x), \zeta(x)$. Given $y\in O_1$ we have the unique expression $\zeta(y)=\omega^{2k'}\phi_3^{\ell'}\phi_4^{m'}(\zeta\sigma^N(y))=\phi_3^{\ell'-1}\phi_4^{k'}\phi_3\phi_4^{m'-k'}(\zeta\sigma^N(y))$ with $0\leq k'\leq n-1, 1\leq \ell'\leq n-1$ and $n\leq m'\leq 2n-1$ depending on $y_{-N},y_{-N+1},\ldots,y_{-1}$. We define $O_1(k',\ell',m')$ as above and then proceed as before. Since $n\geq 3$ there is no solution for $\zeta$ in this case.

Suppose ($ii$), i.e., $C=C\cap D', D\cap D', C'=D\cap C'$ are non-empty. In this case only $\phi_2,\phi_3,\phi_4$ appear when we express $\zeta\sigma^{-1}(x)$ in terms of $\zeta(x)$. Let $x\in O_0$. Then, depending on $x_0,x_1,\ldots,x_{N-1}$, we obtain a unique expression
$$
\zeta\sigma^{-N}(x)=\omega^{j}\phi_2^k\phi_3^\ell\phi_4^m(\zeta(x))=\phi_2^{k-j}\phi_3\phi_2^j\phi_3^{\ell-1}\phi_4^m(\zeta(x))
$$
with $0\leq j\leq n-1, n-1\leq k\leq 2n-2$, and $1\leq\ell,m\leq n$. Then $k+\ell+m=2n$ or $k+\ell+m=3n$ (because $n+1\leq k+\ell+m\leq 4n-2$ and $k+\ell+m=0$ modulo $n$). In the same way as above we define $O_0(j,k,\ell,m)$ as the subset of $O_0$ consisting of $x$ with $x_0,x_1,\ldots,x_{N-1}$ giving this relation on $\zeta\sigma^{-N}(x),\zeta(x)$. Similarly $O_1(j,k,\ell,m)$ is defined as the subset of $O_1$ consisting of $y$ with $y_{-N},y_{-N+1},\ldots,y_{-1}$ giving the relation
$$
\zeta(y)=\omega^j\phi_2^k\phi_3^\ell\phi_4^m(\zeta\sigma^N(y))=\phi_2^{k-j}\phi_3\phi_2^j\phi_3^{\ell-1}\phi_4^m(\zeta\sigma^N(y)).
$$
When $k+\ell+m=2n$ let
$$
L(j,k,\ell,m)=(D\cap D')^{2n-k+j}\times (C\cap D')^{n-1}\times (D\cap D')^{n-j}\times (C\cap D')^{n-\ell+1}\times (D\cap C')^{n-m}
$$
as a subset of $\Lambda^{4n}$, and when $k+\ell+m=3n$ let $L(j,k,\ell,m)=(D\cap D')^{3n-k+j}\times (C\cap D')^{n-1}\times (D\cap D')^{n-j}\times (C\cap D')^{n-\ell+1}\times (D\cap C')^{n-m}$ as a subset of $\Lambda^{4n}$ where the exponent of $D\cap D'$ is increased by $n$. Let $F$ denote the cylinder subset of $X$ determined by
$$
\bigcup O_0(j,k,\ell,m)\times L(j,k,\ell,m)\times L(j',k',\ell',m')\times O_1(j',k',\ell',m')
$$
as a subset of $\prod_{i=-N}^{4N+8n}\Lambda$. Then one shows as before that $K_0\cap \sigma^{2N+8n}(K_1)\cap F$ is not a null set for a sufficiently small $\epsilon>0$ and  $\zeta(\sigma^{-2N-8n}(x))=\zeta(x),\ x\in F$ (e.g.,
$$
\zeta\sigma^{-N-4n}(x)=\phi_4^{n-m}\phi_3^{n-\ell+1}\phi_2^{n-j}\phi_3^{n-1}\phi_2^{2n-k+j}\phi_2^{k-j}\phi_3\phi_2^j\phi_3^{\ell-1}\phi_4^m(\zeta(x))=\zeta(x)
$$
where the first five $\phi$'s are derived from $L(j,k,\ell,m)$). Since this is a contradiction we conclude that $\zeta(x)=\zeta_0$ almost everywhere. This implies that  $\zeta_0=v\zeta_0v^*, \zeta_0=u\zeta_0v^*, \zeta_0=v\zeta_0 u^*$, entailing $uv^*=1$, a contradiction. Thus there is no such $\zeta$.

Suppose ($iii$), i.e., $\zeta\sigma^{-1}(x)$ is one of $\phi_1(\zeta(x)),\phi_2(\zeta(x)),\phi_4(\zeta(x))$. If $x\in O_0$ we deduce that $\zeta\sigma^{-N}(x)$ is uniquely expressed as $\omega^{j}\phi_1^k\phi_2^\ell\phi_4^m(\zeta(x))$ with $0\leq j\leq n-1, 1\leq k\leq n,\, n-1\leq \ell\leq 2n-2$, and $ 1\leq m\leq n$, depending on $x_0,x_1,\ldots,x_{N-1}$. Then $k+\ell+m=2n\ {\rm or}\ 3n$. Note
$$
\zeta\sigma^{-N}(x)=\omega^{j}\phi_1^k\phi_2^\ell\phi_4^m(\zeta(x))=\phi_1^k\phi_2^{\ell-j}\phi_4\phi_2^j\phi_4^{m-1}(\zeta(x)).
$$
Define $O_0(j,k,\ell,m)$ and $O_1(j,k,\ell,m)$ as before and
$$
L(j,k,\ell,m)=(C\cap C')^{n-k}\times (D\cap D')^{2n-\ell+j}\times (D\cap C')^{n-1}\times (D\cap D')^{n-j}\times (D\cap C')^{n-m+1}
$$
as a subset of $\Lambda^{4n}$ when $k+\ell+m=2$ and $L(j,k,\ell,m)$ by the same product as above with the first factor replaced by $(C\cap C')^{2n-k}$ when $k+\ell+m=3n$. We can then proceed as before.

We can treat the cases ($iv$), ($v$), and ($vi$) similarly; so we omit the details.

\begin{prop}\label{Bernoulli-equiv}
Let $X=\prod_{k\in\Z}\Lambda$, $\sigma$, $\mu$, and $u,v\in M_n$ with $n>1$ be as above. Let $C_1,C_1'$ be non-empty proper subsets of $\Lambda$ and let $C,C'$ be the corresponding cylinder subsets of $X$ determined at $0\in\Z$. Define $W=\chi_C\otimes u^*+\chi_D\otimes v^*$ and $W'=\chi_{C'}\otimes u^*+\chi_{D'}\otimes v^*$ in $L^\infty(X,\mu)\otimes M_n$ with $D=X\setminus C$ and $D'=X\setminus C'$ for $C\not=C'$. Then $WV$ is unitarily equivalent to $W'V$ if and only if $n=2$ and $C=D'$ and $D=C'$.
\end{prop}

For $\lambda=(\lambda_1,\lambda_2)\in \T^2$ we define $W(\lambda)=\chi_C\otimes \lambda_1 u^*+\chi_D\otimes \lambda_2 v^*$. Since $\Ad(W(\lambda)V)=\Ad(WV)$ with $W=W(1,1)$ as above, we may ask when $W(\lambda)V$ and $W(\lambda')V$ are unitarily equivalent for $\lambda,\lambda'\in \T^2$.

Suppose that there is a unitary $\zeta\in L^\infty(X,\mu)\otimes M_n$ such that $W(\lambda)V=\zeta W(\lambda')V\zeta^*$. Then we deduce
$$
W(\lambda)(x)=\zeta(x)W(\lambda')(x)\zeta(\sigma^{-1}(x))^*.
$$
Note that $(-1)^{n-1}\det(W(\lambda)(x))=\lambda_1^n\chi_C(x)+\lambda_2^n\chi_D(x)$. Setting $f(x)=\det(\zeta(x))$ and $h(x)=(\lambda'_1\lambda_1^{-1})^n\chi_C(x)+ (\lambda'_2\lambda_2^{-1})^n\chi_D(x)$, as measurable functions on $X$ of modulus one,  we obtain $f(\sigma^{-1}(x))=h(x)f(x)$. Hence for any $N\in\N$
$$
f(x)=\prod_{i=1}^{N} h(\sigma^{i}(x))\cdot f(\sigma^N(x)),\ \ \ f(x)=\prod_{i=0}^{N-1}\overline{h(\sigma^{-i}(x))}\cdot f(\sigma^{-N}(x)).
$$
Since $h(\sigma^k(x))$ depends only on $x_{-k}$  the first equality implies that $f$ is a function measurable with respect to (the Borel sets generated by cylinder sets coming from) $\prod_{i=-\infty}^{-1}\Lambda$ and the second implies that $f$ is a function measurable with respect to $\prod_{i=0}^\infty\Lambda$. (For example if we approximate $f$ by a cylinder function $g$ measurable with respect to $\prod_{i=-N}^{N-1}\Lambda$ in the sense that $\|f-g\|_1<\epsilon$ then it follows that $\|f-\prod_{i=1}^N h\sigma^i\cdot g\sigma^N\|_1<\epsilon$ and $\prod_{i=1}^Nh\sigma^i\cdot g\sigma^N$ is measurable with respect to $\prod_{i=-\infty}^{-1}\Lambda$.) Since $f$ is both measurable with respect to $\prod_{i=-\infty}^{-1}\Lambda$ and $\prod_{i=0}^\infty\Lambda$ we conclude that $f(x)$ is a constant, which implies that $h(x)=1$, i.e., $(\lambda'_1\lambda_1^{-1})^n=1=(\lambda'_2\lambda_2^{-1})^n$.

\begin{prop}\label{Bernoulli-equiv1}
In the situation of Proposition \ref{Bernoulli-equiv} define $W(\lambda)=\chi_C\otimes\lambda_1u^*+\chi_D\otimes \lambda_2v^*$ for $\lambda\in\T^2$ and
let $\lambda,\lambda'\in \T^2$. Then $W(\lambda)V$ is unitarily equivalent to $W(\lambda')V$ if and only if $\lambda_1^n=(\lambda'_1)^n$ and $\lambda_2^n=(\lambda'_2)^n$.
\end{prop}
\begin{pf}
If $\lambda_1^n=(\lambda'_1)^n$ and $\lambda_2^n=(\lambda'_2)^n$ then $\lambda_1=\omega^k\lambda'_1$ and $\lambda_2=\omega^\ell\lambda'_2$ for some $k,\ell\in\Z$. By taking $u^\ell v^{-k}$ for $\zeta$ it follows that $W(\lambda'_1,\lambda'_2)V$ are unitarily equivalent to $W(\omega^{k}\lambda'_1,\omega^{\ell}\lambda'_2)V$. Thus the 'if' part is obvious. The 'only if' part is shown before this proposition.
\end{pf}

\medskip

The other example is based on a dynamical system on $X=\T=\R/\Z$. Let $\theta\in (0,1/2)$ be an irrational number and denote by $\sigma$ the translation by $\theta$: $x\mapsto x+\theta$ on $\T$. If $\alpha$ denotes the automorphism of $C(\T)$ defined by $\alpha(f)(x)=f\sigma^{-1}(x)$, then the C$^*$-algebra crossed product of $C(\T)$ by $\alpha$ is a so-called irrational rotation algebra. Let $\mu$ be the Lebesgue measure on $\T$, which is an ergodic $\sigma$-invariant probability measure and the only $\sigma$-invariant probability measure. Note that $\mu$ is also invariant under the action of $\T$ by translations, which is the fact we will use later. There are many  singular continuous probability measures on $\T$ which are ergodic $\sigma$-quasi-invariant; we will construct such measures in the next section but we do not know if there are ergodic extensions for such measures.

In this case we do not have any results for a C$^*$-version of ergodic extensions since our choice of $W$, similar to the one in the previous case, is not continuous on $\T$.

Let $n$ be an integer greater than $1$ and let $u,v$ be unitaries in $M_n$ as above. Let $C\subset \T$ be a measurable subset of $\T$ such that $0<\mu(C)<1$ and let $D=\T\setminus C$. Define a unitary $W\in M_n$ by $W=\chi_C\otimes u^*+\chi_D\otimes v^*$ and define an automorphism $\beta=\Ad(WV)$ on $L^\infty(\T)\otimes M_n$ (where $V$ is the unitary on $L^2(\T)$ implementing $\alpha$).

\begin{lem}
If $T\in (L^\infty(\T)\otimes M_n)^\beta$ then there is a $T_0\in M_n$ such that $T$ takes values in $\Delta=\{\Ad(u^pv^{-p})(T_0)\ |\ 0\leq p<n\}$ for almost all $x\in\T$.
\end{lem}
\begin{pf}
It follows from the proof of Lemma \ref{finite} that $T$ takes values in $\Gamma=\{\Ad(u^pv^q)(T_0)\ |\ 0\leq p,q<n\}$ for some $T_0\in M_n$. We may suppose that $\{x\in\T\ |\ T(x)=T_0\}$ has positive measure and let $A=\{x\in \T\ |\ T(x)\in \Delta\}$. Then it follows that $T\sigma^{-n}(x)\in \Delta$ for almost all $x\in A$ (by repeating $T(\sigma^{-1}(x))=\Ad\,u(T(x))$ or $\Ad\,v(T(x))$ $n$ times). Since $A$ is $\sigma^n$-invariant and $\sigma^n$ is ergodic it follows that $A$ has full measure.
\end{pf}

Let $T_i,i=0,1,\ldots,m$ be all the distinct elements in the (essential) range of $T$. If $m=0$ then $T_0\in\C1$ (as $T_0=T(\sigma^{-1}(x))=\Ad\,u(T_0)$ for $x\in C$ and $=\Ad\,v(T_0)$ for $x\in D$) and this is what we wanted to prove. Let $F_i=T^{-1}(T_i)$, which is a non-null measurable subset such that $\bigcup_{i}F_i$ has full measure.

Let $(k_r/m_r)$ be the sequence of rational numbers obtained from the continued fraction of $\theta$, which is given as
$$
\begin{pmatrix} k_r \\ m_r \end{pmatrix}=\begin{pmatrix} k_{r-1}& k_{r-2}\\ m_{r-1}& m_{r-2}\end{pmatrix}\begin{pmatrix} b_r\\ 1\end{pmatrix}
$$
for $r\geq 0$ where $(b_r)_{r\geq1}$ is some sequence of natural numbers and $b_0=0, k_{-2}=0,\,k_{-1}=1,\,m_{-2}=1,\,m_{-1}=0$. Note that for $r\geq1$
$$
|\theta-k_r/m_r|<1/m_r^2,\ \ k_r m_{r-1}-k_{r-1}m_{r}=(-1)^{r+1}.
$$
Since $m_r\theta+\Z$ converges to $0\in\T$, the sum $\sum_{i=0}^{m-1}\mu(F_i\cap(F_i+m_r\theta))$ converges to 1 as $r\to \infty$. Let $(r(j))$ be a subsequence such that
$$
\sum_j \big( 1-\sum_i\mu(F_i\cap (F_i+m_{r(j)}\theta)\cap (F_i+m_{r(j)+1}\theta))\big)<\infty
$$ and let $X_0=\bigcup_{\ell=1}^\infty \bigcap_{j=\ell}^\infty G_j$ with $G_j=\bigcup_i (F_i\cap (F_i+m_{r(j)}\theta)\cap (F_i+m_{r(j)+1}\theta))$, which has full measure. Since if $x\in G_j$ then $x,x-m_{r(j)}\theta,x-m_{r(j)+1}\theta\in F_i$ for some $i$ we deduce that
$$
T(\sigma^{-m_{r(j)}}(x))= T(\sigma^{-m_{r(j)+1}}(x))=T(x).
$$
If $x\in X_0$ the equalities hold for all large $j$.

Let $c(x,r)=\#\{i\ |\ x-i\theta\in C,\ 0\leq i<m_r\}$ and $d(x,r)=m_r-c(x,r)$.  Then $ T(\sigma^{-m_{r}}(x))=\Ad(u^{c(x,r)}v^{d(x,r)})(T(x))$ for almost all $x$.

Now we assume that $C$ is an interval $[0,\theta)$ of $\T$. Then it follows that $\mu(\{x\in \T\ |\ c(x,r)=k_r\})= 1-|m_r\theta-k_r|$ (since applying $\sigma$ to $x\in \T$ $m_r$ times results in rotating $x$ around the circle $\T$ almost $k_r$ times; see the lemma below for details). Hence it follows that
$$
T(\sigma^{-m_r}(x))=\Ad(u^{k_r}v^{m_r-k_r})(T(x))
$$
except for $x$ in a subset of measure $|m_r\theta- k_r|$.

Suppose that
$$
\Ad(u^{k_r}v^{m_r-k_r})(T(x))=\Ad(u^{k_{r+1}}v^{m_{r+1}-k_{r+1}})(T(x))=T(x),
$$
which holds for all large $r=r(j)$ for almost all $x\in X_0$.
Since
$$
\begin{pmatrix} k_{r+1}& k_{r}\\ m_{r+1}& m_{r}\end{pmatrix}
$$
has determinant 1 or -1 it follows that there is an inverse matrix consisting of integers, say
$$
\begin{pmatrix} a& b\\ c& d\end{pmatrix}.
$$
(Actually $a=(-1)^{r+1}m_r,\ b=-(-1)^{r+1}k_r,\ c=-(-1)^{r+1}m_{r+1},\ d=(-1)^{r+1}k_{r+1}$.)
Since $(u^{k_r}v^{m_r-k_r})^c(u^{k_{r+1}}v^{m_{r+1}-k_{r+1}})^a$ is proportional to $uv^*$ (as the exponent of $u$ is $k_r c+ k_{r+1}a=1$ and the exponent of $v$ is $m_rc-k_rc+m_{r+1}a-k_{r+1}a=-1$), it follows that $\Ad(uv^*)(T(x))=T(x)$. Since $T(x)=\Ad(u^pv^{-p})(T_0)$ for some $p$, this implies that $\Ad(uv^*)(T_0)=T_0$, i.e., $m=0$ and $T(x)=T_0$ for almost all $x$. Since $T\sigma^{-1}(x)=\Ad\,u(T(x))$ or $\Ad\,v(T(x))$ depending on $x$, this shows that $T_0\in\C1$. Thus we conclude that $(L^\infty(\T)\otimes M_n)^\beta=\C1$.

\begin{lem}
Let $m\in\N$ and let $[m\theta]$ be the largest integer satisfying $[m\theta]\leq m\theta$. Then
\begin{align*}
&\mu(\{x\in \T\ |\ c'(x,m)=[m\theta]\})=[m\theta]+1-m\theta,\\
&\mu(\{x\in\T\ |\ c'(x,m)=[m\theta]+1\})= m\theta-[m\theta],
\end{align*}
 where $c'(x,m)=\#\{i\ |\ x-i\theta\in C, 0\leq i<m\}$. If $m\theta\approx k$ then $\mu(\{x\in \T\ |\ c'(x,m)=k\})=1-|m\theta-k|$.
\end{lem}
\begin{pf}
We have assumed that $C=[0,\theta)$. Since $C$ is an interval of length $\theta<1/2$, if $x-i\theta\in C$ then $x-(i-1)\theta\not\in C$ and $x-(i+1)\theta\not\in C$ and if two consecutive points $x-i\theta,x-(i+1)\theta$ in orbit passes the middle point $\theta/2$ of $C$ then one and only one of them falls into $C$. If $m'$ is the smallest positive integer with $[m'\theta]=[m\theta]$ then $c'(x,m')\geq [m\theta]$ (as at least $[m\theta]$ of $x,x-\theta,\ldots,x-(m'-1)\theta$ belong to $C$) and $c'(x,m')=[m\theta]+1$ occurs only when $(m'-1)\theta,x\in C$. If $c'(x,m')=[m\theta]$ and one of $x-m'\theta,x-(m'+1)\theta\ldots,x-(m-1)\theta$ belongs to $C$ then $c'(x,m)=[m\theta]+1$.

If $[(m-1)\theta]=[m\theta]$ then $m'\leq m-1$. Note that $x$ lies on the arc of length $\theta$ between $(m'-1)\theta$ and $m'\theta$. Then it follows that if one of $x,x-m'\theta,x-(m'+1)\theta,\ldots,x-(m-1)\theta$  belongs to $C$ then $c'(x,m)=[m\theta]+1$; otherwise $c'(x,m)=[m\theta]$. The set of $x$ with $c(x,m)=[m\theta]+1$ is given by the condition $0\leq x<\theta$, $m'\theta-[m\theta]\leq x<(m'+1)\theta-[m\theta], \ldots,\ {\rm or}\ (m-1)\theta-[m\theta]\leq x<m\theta-[m\theta]$, i.e., $0\leq x<m\theta-[m\theta]$.

If $[(m-1)\theta]<[m\theta]$ (i.e., $[(m-1)\theta]<(m-1)\theta<[m\theta]<m\theta$) then $m'=m$ and only when both $x,x-(m-1)\theta$ fall into $C$ we have that $c'(x,m)=[m\theta]+1$. Since $0\leq x<m\theta-[m\theta]$ if and only if $0<[m\theta]-(m-1)\theta\leq x+[m\theta]-(m-1)\theta<\theta$, the set of $x$ with $c(x,m)=[m\theta]+1$ is $0\leq x<m\theta-[m\theta]$.

If $k>m\theta>[m\theta]$ then $k=[m\theta]+1$; so $m\theta-[m\theta]=1-(k-m\theta)$ and if $m\theta>k$ then $k=[m\theta]$; $[m\theta]+1-m\theta=1-(m\theta-k)$. Thus we derive the last statement.
\end{pf}

We have shown the following:

\begin{prop}\label{irrational}
Let $\sigma$ denote the map $x\mapsto x+\theta$ on $\T=\R/\Z$ where $\theta\in (0,1/2)$ is irrational and let $\mu$ be the Lesbegue measure on $\T$. Let $V$ denote the the unitary on $L^2(\T)$ induced by $\sigma$. Let $n$ be an integer greater than 1 and let $W=\chi_C\otimes u^*+\chi_D\otimes v^*$ where $C$ is an interval $[0,\theta)\subset \T$ and $D=\T\setminus C$ where $u,v$ are the canonical pair of unitaries in $M_n$ defined before. Then the automorphism $\Ad(WV)$ acts on $L^\infty(\T)\otimes M_n$ ergodically.
\end{prop}

In the case $n=1$ we may choose a scaler for $W$ and ask when $\lambda V $ is unitarily equivalent to $\lambda'V$ for $\lambda,\lambda'\in\T=\{z\in\C\ |\ |z|=1\}$. That is, when is there a unitary $\zeta\in L^\infty(\T)$ such that $\lambda V=\zeta \lambda'V\zeta^*$, or $\lambda\overline{\lambda'}=\zeta(x)\overline{\zeta(x-\theta)},\ x\in\T$? Let $\lambda\overline{\lambda'}=e^{2\pi i\eta}$ with $\eta\in\R$. Since $e^{2\pi i q\eta}=\zeta(x)\overline{\zeta(x-q\theta)}$ for $q\in\Z$ and $\zeta(\,\cdot\,)\overline{\zeta(\,\cdot\,-q\theta)}$ converges to 1 as $q\theta$ converges to 0 in $\T$, we deduce the hypothesis  of the following lemma.

\begin{lem}\label{irrationalnumbers}
Let $\theta,\eta\in \R$. Suppose that $\theta\in (0,1/2)$ is an irrational number and that if $q_k\theta$ converges to $0$ in $\T=\R/\Z$ then
$q_k\eta$ converges to $0$ in $\T$ for any sequence $(q_k)$ in $\Z$. Then $\eta=m\theta$ for some $m\in\Z$.
\end{lem}
\begin{pf}
Define a map $\phi$ of $\theta \Z/\Z\subset \T=\R/\Z$ into $\T$ by $q\theta+\Z\mapsto q\eta+\Z,\ q\in\Z$. This is well-defined because $\theta$ is irrational. If $(q_k)$ is a sequence in $\Z$ and $(q_k\theta+\Z)$ is Cauchy then $(q_k\eta+\Z)$ is Cauchy too. (If $(q_k\eta+\Z)$ is not Cauchy then there are subsequences $k(\ell), k'(\ell)$ such that $(q_{k(\ell)}-q_{k'(\ell)})\eta+\Z$ does not converges to zero, which contradicts that $(q_{k(\ell)}-q_{k'(\ell)})\theta+\Z$ converges to zero.) Hence $\phi$ extends to a continuous map of $\T$ into $\T$. Since $\phi(\T)$ is a connected compact subset of $\T$, $\eta$ is either $0$ or an irrational. If $\mu$ is irrational then $\phi$ is onto. If $x=q\theta$ and $y=q'\theta$ with $q,q'\in\Z$, then it follows that $\phi(x+y)=q\eta+q'\eta,\ {\rm mod}\ \Z$. Hence we deduce that $\phi(x+y)=\phi(x)+\phi(y)$ in $\T$ for all $x,y\in \T$ and that $\phi(rx)=r\phi(x), x\in\T$ for all rational $r$. Since $\phi$ is continuous we have $\phi(tx)=t\phi(x)$ for all $t\in\R$. Since $\phi(0)=0$ there is a neighborhood $U$ of $0\in\T$ (identified with $(-1/2,1/2]$) such that $U\subset (-1/2,1/2)$ and  $\phi(U)\subset (-1/2,1/2)$. Let $x\in U\setminus\{0\}$. Then $\phi(tx)=t\phi(x)=mtx$ with $m=\phi(x)/x$, or $\phi(t)=mt$. Since $\phi$ is a map from $\T$ onto $\T$, it follows that $m$ is an integer.
\end{pf}

\begin{prop}\label{equiv0}
In the situation of Proposition \ref{irrational} suppose that $n=1$ $($and $W=1${\rm)}.
The following conditions are equivalent for $\lambda,\lambda'\in\T$.
\begin{enumerate}
\item $\lambda V$ and $\lambda'V$ are unitarily equivalent.
\item $\lambda=e^{2\pi im\theta}\lambda'$ for some $m\in\Z$.
\end{enumerate}
\end{prop}
\begin{pf}
Define $Y\in C(\T)\subset L^\infty(\T)$ by $Y(x)=e^{2\pi i x}$. Then $YVY^*=e^{2\pi i\theta}V$. Hence $(2)\Rightarrow$(1). The other implication follows from Lemma \ref{irrationalnumbers} and its preceding remark.
\end{pf}

The above proposition for $n=1$ is perhaps known. We present a version for $n>1$ in the situation of Proposition \ref{irrational}.

Let $\lambda=(\lambda_1,\lambda_2)\in\T^2$  and let
$$
W(\lambda)=\chi_C\otimes \lambda_1u^*+\chi_D\otimes \lambda_2 v^*.
$$
Then it follows that $\Ad(W(\lambda)V)=\Ad(W(\lambda')V)$ on $L^\infty(\T)\otimes M_n$ for all $\lambda'\in \T^2$ (as $W(\lambda)=(\chi_C\otimes \lambda_11+\chi_D\otimes\lambda_21)W$ and $\chi_C\otimes\lambda_1 1+\chi_D\otimes \lambda_2 1$ is in $L^\infty(\T)\otimes \C1$); so we could ask when $W(\lambda)V$ and $W(\lambda')V$ are unitarily equivalent, i.e.,  there is a unitary $\zeta\in L^\infty(\T)\otimes M_n$ such that $W(\lambda)V=\zeta W(\lambda')V\zeta^*$. By taking $\zeta=1\otimes u^kv^\ell$ it follows that $W(\lambda)V$ is unitarily equivalent to $W(\omega^{\ell}\lambda_1,\omega^{k}\lambda_2)$ for all $k,\ell$ with $\omega=e^{2\pi i/n}$. Similarly by taking $\zeta=U_a\otimes 1$ with $U_a(x)=e^{2\pi iax},\ x\in\T$ for some $a\in\R$ it follows that $W(\lambda)V$ is unitarily equivalent to $W(e^{2\pi ia(1-\theta)}\lambda_1,e^{2\pi ia\theta}\lambda_2)V$.

\begin{prop}\label{equiv1}
In the situation of Proposition \ref{irrational} the following conditions are equivalent for $\lambda,\lambda'\in \T^2$.
\begin{enumerate}
\item $W(\lambda)V$ and $W(\lambda')V$ are unitarily equivalent, i.e., there is a unitary $\zeta\in L^\infty(\T)\otimes M_n$ such that $W(\lambda)=\zeta W(\lambda')\bar{\alpha}(\zeta)^*$ where $\bar{\alpha}=\Ad\,V$.
\item $\lambda_1^n=e^{2\pi ia(\theta-1)}(\lambda'_1)^n$ and $\lambda_2^n=e^{2\pi i a\theta}(\lambda'_2)^n$ for some $a\in\R$.
\end{enumerate}
\end{prop}
\begin{pf}
Suppose (2). Then $\lambda_1=e^{2\pi ia(\theta-1)/n}\omega^k\lambda'_1$ and $\lambda_2=e^{2\pi i a\theta/n}\omega^\ell\lambda_2'$ for some $k,\ell\in\Z$. Since $W(\lambda_1',\lambda'_2)V$ is unitarily equivalent to $W(\omega^k\lambda'_1,\omega^\ell\lambda'_2)V$ we may suppose, replacing $a/n$ by $a$ too, that $\lambda_1=e^{2\pi ia(\theta-1)}\lambda'_1$ and $\lambda_2=e^{2\pi ia\theta}\lambda'_2$. We have shown in this case $W(\lambda)V$ and $W(\lambda')V$ are unitarily equivalent just before this proposition.

Suppose (1). Then it follows that $W(\lambda)(x)=\zeta(x)W(\lambda')(x)\zeta(x-\theta)^*$ for almost all $x\in \T$. Since $\det(W(\lambda)(x))=\lambda_1^n(-1)^{n-1}\chi_C(x)+\lambda_2^n (-1)^{n-1}\chi_D$ (as $\det(u)=\det(v)=(-1)^{n-1}$) we deduce that
\begin{align*}
\lambda_1^n=(\lambda'_1)^n\det(\zeta(x))\overline{\det(\zeta(x-\theta))}, &\ \  x\in C,\\
\lambda_2^n=(\lambda_2')^n\det(\zeta(x))\overline{\det(\zeta(x-\theta))}, &\ \  x\in D.
\end{align*}
Define $\eta_i\in (-\pi,\pi]$ by $e^{2\pi i\eta_i}=(\lambda_i\overline{\lambda'_i})^n$ for $i=1,2$.
Setting $f(x)=\det(\zeta(x))$ and $h(x)= e^{2\pi i\eta_1}\chi_C(x) + e^{2\pi i\eta_2}\chi_D(x)$ as measurable functions of modulus one on $\T$, this amounts to
$$f(x)\overline{f(x-\theta)}=h(x).
$$
If $q_n\theta-p_n$ converges to zero with $q_n,p_n\in\Z$ then $f(x)\overline{f(x-q_n\theta)}\to 1$ in $L^1$. Since
$$
f(x)\overline{f(x-q_n\theta)}=\prod_{k=0}^{q_n-1}h(x-k\theta),
$$
whose right-hand side is $e^{2\pi i(p_n\eta_1+(q_n-p_n)\eta_2)}$ outside a subset of small measure, one concludes that $p_n\eta_1+(q_n-p_n)\eta_2$ converges to zero in $\T$.

As in the proof of Lemma \ref{irrationalnumbers} there is a continuous map $\phi$ of $\R$ into $\T$ such that $\phi(q\theta-p)=p\eta_1+(q-p)\eta_2=\eta_2q+(\eta_1-\eta_2)p$. Since $\phi$ satisfies that $\phi(x+y)=\phi(x)+\phi(y),\ x,y\in\R$, we conclude that $\phi(x)=ax$ (mod $\Z$) where $a$ is a constant. If $a=0$ then $\eta_i=0$. Suppose that $a\not=0$. Since $a(q\theta -p)=\eta_2 q+(\eta_1-\eta_2)p$ (mod $\Z$) we obtain that $a\theta=\eta_2$ and $a=\eta_2-\eta_1$ or $\eta_1=-a(1-\theta)$ and $\eta_2=a\theta$ (mod $\Z$). This concludes the proof.
\end{pf}

If we restrict ourselves to the case $\lambda_1=\lambda_2$ in the above proposition, then it follows that $a$ must be an integer, i.e., $\lambda WV$ is unitarily equivalent to $\lambda' WV$ if and only if $\lambda^n=e^{2\pi i m\theta}(\lambda')^n$ for some $m\in\Z$ (cf. Proposition \ref{equiv0}).

\medskip
We have to leave many problems unanswered. For example we did not explore all possible ergodic extensions in Propositions \ref{Bernoulli} and \ref{irrational} in the case of $M_n=\B(\C^n)$, let alone the case of $\B(\Hil)$ with $\Hil$ an infinite-dimensional Hilbert space. We did not attempt to solve the problem for general  $(X,\sigma)$.


\section{Quasi-invariant measures}

Let $\theta\in (0,1)$ be an irrational number and let $\sigma$ denote the homeomorphism on $\T=\R/\Z$ defined by $x\mapsto x+\theta$ (mod $1$). In this case we have noted that there are at least two kinds of ergodic $\sigma$-quasi-invariant probability measures on $\T$: the Lebesgue measure on $\T$ (which is $\sigma$-invariant) and an atomic measure on each orbit $\{x+m\theta\ |\ m\in\Z\}$.

We shall construct ergodic $\sigma$-quasi-invariant singular continuous probability measures on $\T$. In the following we denote by $(x)$ the representative in $(-1/2,1/2]$ of $x+\Z\in\T$ for $x\in\R$.

Let $P=\prod_{i=1}^\infty\{0,1\}$, the infinite direct product of copies of $\{0,1\}$ with the product topology. We define a continuous map $\Phi:P\to \T$ as follows: Let $(m_i)$ be an increasing sequence in $\N$ such that $|(m_1\theta)|<1/3$ and $|(m_i\theta)|<|(m_{i-1}\theta)|/3$ for $i>1$. With such a sequence $(m_i)$ let
$$
\Phi(x)=\sum_{i=1}^\infty x_i(m_i\theta),\ \ x=(x_i)\in P.
$$
Then $\Phi$ is well-defined and continuous.

\begin{lem}
Let $\lambda_i\in \{-1,0,1\}$ for $i=1,2,\ldots$. If $t=\sum_i\lambda_i(m_i\theta)$ then $|t|<1/2$ and $(\lambda_i)$ is uniquely determined by $t$.
\end{lem}
\begin{pf}
If $t=\sum_i\lambda_i(m_i\theta)$ then $|t|\leq \sum_i|(m_i\theta)|<1/2$.

If $t=\sum_i\lambda_i'(m_i\theta)$ for another sequence $(\lambda'_i)$ in $\{-1,0,1\}$ different from $(\lambda_i)$ then there is $N\in \N$ such that $\lambda_i=\lambda'_i$ for $i<N$ and $\lambda_N\not=\lambda'_N$. Since $\sum_{i=N}^\infty (\lambda_i-\lambda'_i)(m_i\theta)=0$ it follows that
$|(m_N\theta)|\leq 2\sum_{i=N+1}^\infty |(m_i\theta)|.$
But from the assumption on $(m_i)$ we deduce that
$$
\sum_{i=N+1}^\infty |(m_i\theta)|<\sum_{i=1}^\infty |(m_N\theta)|/3^i=|(m_N\theta)|/2,
$$
which is a contradiction.
\end{pf}

The above lemma, in particular, implies that  $\Phi$ is injective. Hence we conclude that $\Phi(P)$ is a compact subset of $\T$ and $\Phi$ is a homeomorphism of $P$ onto $\Phi(P)$. Let $a_N$ denote the sum of $(m_i\theta)<0$ with $i>N$ and $b_N$ the sum of $(m_i\theta)>0$ with $i>N$. (It follows from the above calculation that $b_N-a_N=\sum_{i=N+1}^\infty |(m_i\theta)|<|(m_N\theta)|/2$.) Since $\Phi(P)$ is contained in the intersection of the decreasing sequence
$$
\sum_{S\subset \{1,2,\ldots,N\}} \sum_{i\in S}(m_i\theta)+[a_N,b_N]
$$
of closed subsets of $\T$ (as each closed interval $\sum_{i\in S}(m_i\theta)+[a_N,b_N]$ shrinks into two disjoint intervals $\sum_{i\in S}(m_i\theta)+[a_{N+1},b_{N+1}]$ and its translate by $(m_{N+1}\theta)$ at the next stage and $b_N-a_N=|(m_{N+1}\theta)|+ b_{N+1}-a_{N+1}>3(b_{N+1}-a_{N+1})$), $\Phi(P)$ is a Cantor set. Thus it follows that $\Phi(P)$ is a null set with respect to the Lebesgue measure.

Let $\nu_0$ denote the product measure on $P$ given by $\prod_{i=1}^\infty \{a_i,1-a_i\}$ for some sequence $(a_i)$ with $0<a_i<1$ such that $\sup a_i<1$ and $\inf a_i>0$. Note that $\nu_0$ is a non-atomic measure (even when restricted to the subfield generated by the co-ordinates $x_{n_1}, x_{n_2},\ldots$ for any subsequence $(n_i)$). We define a probability measure $\nu_0'$ on $\T$ by $\nu_0': A\mapsto \nu_0\Phi^{-1}(A\cap \Phi(P))$ and then a probability measure $\nu$ on $\T$ as follows:
$$
\nu(A)=\sum_{k=-\infty}^\infty \gamma^{1+|k|} \nu_0'\sigma^k(A),
$$
where $\gamma=\sqrt{2}-1$. Let $B=\bigcup_k \sigma^k(\Phi(P))$, an $F_\sigma$ subset of $\T$. Then $\nu(B)=1$ and $B$ has Lesbegue measure $0$. Since $\nu_0'$ is non-atomic, so is $\nu$. Thus $\nu$ is a non-atomic measure singular from the Lebesgue measure, i.e., $\nu$ is {\em singular continuous}. Since $\gamma\nu(A)\leq \nu\sigma(A)\leq \gamma^{-1}\nu(A)$ for all Borel sets $A$ we conclude that $\nu$ is $\sigma$-quasi-invariant. Note that $\nu$ is not $\sigma$-invariant and is not equivalent to a $\sigma$-invariant probability measure. (If it is $\sigma$-invariant and $\phi$ is the state on $C(\T)$ defined by $\nu$ and $Y\in C(\T)$ is defined by $Y(t)=e^{2\pi it}$, then one can show that $\phi(Y^k)=\phi( N^{-1}\sum_{k=0}^{N-1}\alpha^k(Y))$, which is valid for all $N$ yielding $\phi(Y^k)=0$ for $k\not=0$ as $\alpha(Y)=Y\sigma^{-1}= e^{-2\pi i\theta}Y$. Thus it follows that $\phi(f)=\int fdt$ for $f\in C(\T)$, a contradiction.)

\begin{lem}\label{int}
Let $t\in \R$. Then $\Phi(P)\cap (\Phi(P)+t)$ is a null set with respect to $\nu'_0$ if and only if $t$ does not belong to $\{\sum_i\lambda_i (m_i\theta)\ |\ \lambda_i=-1,0,+1\}$ modulo $\Z$ $($where $\lambda_i=0$ except for a finite number of $i)$.
\end{lem}
\begin{pf}
If $\Phi(P)\cap (\Phi(P)+t)\not=\emptyset$ then there are $x,y\in P$ such that $\Phi(x)=\Phi(y)+t$, or $t=\sum_i(x_i-y_i)(m_i\theta)+k$ for some $k\in\Z$.  From the previous lemma it follows that $t-k$ determines $\lambda_i=x_i-y_i\in \{-1,0,1\}$. If $\lambda_i=0$ except for a finite number of $i$ then $t=\sum_i\lambda_i m_i\theta$ modulo $\Z$ and $\Phi(P)\cap (\Phi(P)+t)$ is equal to the set of $\Phi(x)$ satisfying $x_i=1$ if $\lambda_i=1$, $x_i=0$ if $\lambda_i=-1$, and $x_i$ is arbitrary if $\lambda_i=0$, which implies that $\Phi(P)\cap (\Phi(P)+t)$ has positive measure. Otherwise $\Phi(P)\cap(\Phi(P)+t)$ is a null set.
\end{pf}

\begin{lem}\label{nu0}
$\nu|\sigma^k(\Phi(P))$ is mutually absolutely continuous with respect to $\nu_0'|\sigma^k(\Phi(P))$ for all $k\in\Z$.
\end{lem}
\begin{pf}
We may suppose that $k=1$ (as $\nu$ is a kind of average of $\nu_0\sigma^k$ over $k$).
It is obvious that $\nu'_0|\Phi(P)$ is absolutely continuous with respect to $\nu|\Phi(P)$.  What we have to show is that $\nu'_0\sigma^k|\Phi(P)$ is absolutely continuous with respect to $\nu'_0|\Phi(P)$ for any $k$. It follows from the previous lemma that $k$ must be $\sum_i\lambda_im_i$ for some $(\lambda_i)$ in order that $\nu'_0\sigma^k|\Phi(P)$ is non-zero.

In this case $\Phi(P)\cap \sigma^k(\Phi(P))$ is the cylinder set of $P$ determined by $x_i=1$ for $i\in I_1=\{i\ |\ \lambda_i=1\}$ and $x_i=0$ for $i\in I_0=\{i\ |\ \lambda_i=-1\}$. (Here and henceforth we identify $\Phi(P)$ with $P$ and $\nu_0'$ with $\nu_0$.) The inverse image of $\Phi(P)\cap \sigma^k(\Phi(P))$ under $\sigma^k$ is the cylinder set $Q$ of $P$ determined by $x_i=0$ for $i\in I_1$ and $x_i=1$ for $i\in I_0$. Note that $\nu_0(Q)=\prod_{i\in I_1}a_i\cdot \prod_{i\in I_0}(1-a_i)$ and $\nu_0\sigma^k(Q)=\prod_{i\in I_1}(1-a_i)\cdot \prod_{i\in I_0}a_i$. Then the definition of $\nu_0$ implies $\nu_0\sigma^k|Q=c_k\nu_0|Q$ where
$$
c_k=\prod_{i\in I_1} (1-a_i)/a_i \cdot \prod_{i\in I_0} a_i/(1-a_i).
$$
Hence we conclude that $\nu_0\sigma^k|\Phi(P)\leq c_k\nu_0|\Phi(P)$.
\end{pf}

\begin{prop}\label{sing-cont}
Let $\nu$ be a probability measure on $\T$ constructed from $\nu_0$ on $P$ and $\Phi:P\to \T$ as above. Then $\nu$ is an ergodic $\sigma$-quasi-invariant singular continuous measure on $\T$.
\end{prop}
\begin{pf}
What remains to show is that $\nu$ is ergodic with respect to $\sigma$.

Suppose that $A$ is a $\sigma$-invariant measurable subset of $\T$ with $0<\nu(A)$. Let $A_0=A\cap \Phi(P)$, which has positive measure since $\sigma^k(A_0)=A\cap \sigma^k(\Phi(P))$ and $A=\bigcup_k A\cap \sigma^k(\Phi(P))$ (modulo null sets). We regard $A_0$ as a measurable subset of $P$. If $x\in A_0$ and $y\in P$ satisfies $x_i=y_i$ for all large $i$ then $y\in A_0$ (by Lemma \ref{int}), i.e., $A_0$, as a measurable subset of $P$, does not depend on the first $N$-coordinates for any $N\in\N$. This implies that $\nu_0(A_0\cap C)=\nu_0(A_0)\nu_0(C)$ for any cylinder set $C$ of $P$ and hence for any measurable set $C$. Thus we conclude that $\nu_0(A_0)=\nu_0(A_0)^2$, i.e., $\nu_0(A_0)=1$ or $\nu(\Phi(P)\setminus A)=0$, which implies that $\nu(A)=1$. Hence $\nu$ is ergodic.
\end{pf}

Let $a_i=a$ for all $i$ and denote by $\nu_{0a}$ the corresponding probability measure $\nu_0$ on $P$. If $a,b\in (0,1)$ are different then $\nu_{0a}$ and $\nu_{0b}$ are mutually singular. By using the same $\Phi:P\to \T$ we construct a probability measure $\nu_a$ on $\T$ from $\nu_{0a}$. They are all ergodic $\sigma$-quasi-invariant singular continuous probability measures on $\T$.

\begin{cor}
The above $\nu_a, 0<a<1$ are mutually singular.
\end{cor}
\begin{pf}
Let $a,b\in (0,1)$ with $a\not=b$. Since $\nu_a|P\simeq \nu_{0a}|P$ and $\nu_b|P\simeq \nu_{0b}|P$ and $\nu_{0a}$ and $\nu_{0b}$ are mutually singular we deduce that $\nu_a|P$ and $\nu_b|P$ are mutually singular; in particular $\nu_a|P\not=\nu_b|P$. Since $\nu_a$ and $\nu_b$ are both ergodic, we conclude that $\nu_a$ and $\nu_b$ are mutually singular.
\end{pf}

Let $\Sigma=(\T,\sigma)$ and let $\nu$ be an ergodic $\sigma$-quasi-invariant probability measure on $\T$. Then $\pi_{(\nu,\C,1)}$ is a simplest kind of irreducible representations of $C^*(\Sigma)$.

Let $V$ denote the unitary on $L^2(\T,\nu)$ defined by $(V\xi)(x)=\xi\sigma^{-1}(x)(d\nu\sigma^{-1}/d\nu)(x)^{1/2}$.
Then the spectrum of $V$ is $\T$ and there is a probability measure $\nu_1$ on $\T$ such that the isomorphism  $Y^k\mapsto V^k,\ k\in \Z$ of $C(\T)$ into $\B(L^2(\T,\nu))$ extends to the one of $L^\infty(\T,\nu_1)$ into $\B(L^2(\T,\nu))$, where $Y:x\mapsto e^{2\pi ix}$. Since $VYV^*=e^{-2\pi i\theta}Y$ or $YVY^*=e^{2\pi i\theta}V$, we deduce that $\nu_1$ is quasi-invariant under $\sigma$. Let $V_1$ denote the unitary on $L^2(\T,\nu_1)$ defined by $(V_1\xi)(x)=\xi\sigma^{}(x)(d\nu_1\sigma^{}/d\nu)(x)^{1/2}$ (where we have used $\sigma$ instead of $\sigma^{-1}$). Then by Proposition \ref{irreducible} we conclude that $L^2(\T,\nu)\cong L^2(\T,\nu_1)\otimes \Hil$ for some Hilbert space $\Hil$ where $V$ (resp. $Y$) corresponds to $Y\otimes 1$ (resp. $W(V_1\otimes 1)$) for some unitary $W$ in $L^\infty(\T,\nu_1)\otimes \B(\Hil)$. That is, exchanging the roles of $Y$ and $V$ we deduce that $\pi_{(\nu,\C,1)}$ is equivalent to $\pi_{(\nu_1,\Hil,W)}$, an irreducible representation for $(\T,\sigma^{-1})$.

Suppose that $\nu$ is the Lebesgue measure; in this case $V$ has a complete set of eigenvectors. Then $\nu_1$ must be atomic and ergodic. Then by Proposition \ref{cocyclevanishing} we obtain $\Hil\cong\C$ and can assume that $W=1$. The converse also follows.

Suppose that $\nu$ is atomic, i.e., $Y$ has a complete set of eigenvectors. If $\dim(\Hil)>1$ then $W(V_1\otimes 1)$ has no eigenvalues by Proposition \ref{noeigenvectors}, which contradicts that $Y$ is diagonal. Thus $\Hil=\C$ and hence $L^2(\T,\nu)\cong L^2(\T,\nu_1)$. Hence $WV_1$ has an eigenvector, say $WV_1\xi=\lambda\xi$ for a unit vector $\xi\in L^2(\T,\nu_1)$ and a complex number $\lambda$ of modulus 1. Then it follows that
$$
W(x)\xi\sigma(x)(\frac{d\nu_1\sigma}{d\nu_1}(x))^{1/2}=\lambda\xi(x).
$$
Hence we deduce that $|\xi(x)|^2d\nu_1(x)$ is a $\sigma$-invariant probability measure, which must be the Lebesgue measure on $\T$. (In this case $V_1$ is diagonal and hence $W$ must be a constant.) Thus we have:

\begin{prop}
In the above situation  $\nu$ is mutually absolutely continuous with respect to the Lebesgue measure if and only if $\nu_1$ is an ergodic $\sigma$-quasi-invariant atomic probability measure.  Furthermore $\nu$ is atomic if and only if $\nu_1$ is mutually absolutely continuous with respect to the Lebesgue measure and $\Hil=\C$. In these cases $L^2(\T,\nu)\cong L^2(\T,\nu_1)$.
\end{prop}

We do not know if the case $\dim(\Hil)>1$ can actually occur when we start from $\pi=\pi_{(\nu,\C,1)}$ or if $V''$ can fail to be maximal abelian (when $\pi(Y)''$ is maximal abelian). If $\nu$ is singular continuous then $\nu_1$ is either singular continuous or mutually absolutely continuous with respect to the Lebesgue measure with $\dim(\Hil)>1$.


\end{document}